\documentclass[a4paper,11pt]{article}
\usepackage[nointlimits]{amsmath}
\usepackage{amsfonts}
\usepackage{amssymb}
\usepackage{amsmath}
\usepackage{amsthm}
\usepackage{latexsym}
\usepackage{graphicx}
\usepackage{layout}
\usepackage[english]{babel}
\usepackage{color}

\usepackage{esint}

\usepackage{fullpage}

\numberwithin{equation}{section} 
        
\newtheorem{lemma1}     {Lemma}[section]
\newtheorem{teorema1}   [lemma1]{Theorem}
\newtheorem{prop1}      [lemma1]{Proposition}
\newtheorem{coroll1}    [lemma1]{Corollary}
\newtheorem{cong1}      [lemma1]{Conjecture}
\newtheorem{remark1}    [lemma1]{Remark}
\newtheorem{defin1}     [lemma1]{Definition}

\newenvironment{Lemma}[1][]
        {\begin{lemma1}[#1]\begin{samepage}}{\end{samepage}\end{lemma1}}
\newenvironment{Theorem}[1][]
        {\begin{teorema1}[#1]\begin{samepage}}{\end{samepage}\end{teorema1}}
\newenvironment{Proposition}[1][]
        {\begin{prop1}[#1]\begin{samepage}}{\end{samepage}\end{prop1}}
\newenvironment{Corollary}[1][]
        {\begin{coroll1}[#1]\begin{samepage}}{\end{samepage}\end{coroll1}}

\providecommand{\R}{\mathbb{R}}
\providecommand{\N}{\mathbb{N}}

\providecommand{\ud}{\, \mathrm{d}}
\providecommand{\dx}{\, \mathrm{d} x}
\providecommand{\dy}{\, \mathrm{d} y}
\providecommand{\dz}{\, \mathrm{d} z}
\providecommand{\dxy}{\, \mathrm{d} \x}
\providecommand{\dxyz}{\, \dxy}
\renewcommand{\div}{\mathrm{div}}
\newcommand{\eps}{\varepsilon}
\def\lt{\left}
\def\rt{\right}
\newcommand{\F}{F_\ell}
\DeclareMathOperator{\supp}{supp}

\def\Om{\Omega}
\def\les{\lesssim}
\def\ges{\gtrsim}
\def\Omh{\Omega_h}
\def\x{\mathbf{x}}
\renewcommand{\c}{c_0}
\def\Xint#1{\mathchoice
{\XXint\displaystyle\textstyle{#1}}%
{\XXint\textstyle\scriptstyle{#1}}%
{\XXint\scriptstyle\scriptscriptstyle{#1}}%
{\XXint\scriptscriptstyle\scriptscriptstyle{#1}}%
\!\int}
\def\XXint#1#2#3{{\setbox0=\hbox{$#1{#2#3}{\int}$ }
\vcenter{\hbox{$#2#3$ }}\kern-.57\wd0}}

\def\dashint{\Xint-}

\title{Study of island formation in epitaxially strained films on unbounded domains}
\author{P. Bella
\footnote{Max Planck Institute for Mathematics in the Sciences, Leipzig (Germany), email: bella@mis.mpg.de}
\and
        M. Goldman
        \footnote{Max Planck Institute for Mathematics in the Sciences, Leipzig (Germany), email: goldman@mis.mpg.de}
        \and B. Zwicknagl
        \footnote{Institute for Applied Mathematics, Bonn University, Bonn (Germany), email: zwicknagl@iam.uni-bonn.de}}
\begin{document}
\bibliographystyle{plain}

\maketitle
\begin{abstract}
We consider a variational model {{related to}} the formation of islands in heteroepitaxial growth on unbounded domains. We first derive the scaling regimes of the minimal energy
 in terms of the volume of the film and the amplitude of the crystallographic misfit. For small volumes, non-existence of minimizers is then  proven. This corresponds
 to the experimentally observed wetting effect. On the other hand, we show the existence of minimizers for large volumes. We finally study the asymptotic behavior of the optimal shapes.
\end{abstract}

\section{Introduction} 

We consider the epitaxial deposition of a thin crystalline film on a relatively thick rigid substrate with a misfit between the lattice parameters of the film and those of the substrate. Experimental and numerical observations suggest that the shape of the film changes with increasing volume (see \cite{BC:02,GaoNix:99,GrHulFlo:04, SpeTer:10,Smere:12}). At small volumes, one typically observes a very thin flat layer (``wetting''), while at larger volumes, compact islands form. This transition is often explained as the result of a competition between two opposing types of energies, namely, the stored strain energy due to the crystallographic misfit, and the surface energy of the film's free surface. Heuristically, at small volumes, the surface energy dominates, and complex structures are avoided, while at larger volumes, the film forms patterns to release elastic energy at the price of an additional surface energy.

We study analytically a two-dimensional variational model  introduced in \cite{Spen:99} (see also \cite{BC:02,FFLM:07,FM:12}), to describe the surface morphologies of the epitaxially strained film. 
The main difference to the previous analytical works (see \cite{FFLM:07,FM:12,GolZwick}) is that the model explicitly allows for wetting, which corresponds to film profiles with unbounded support.  
We assume that the film occupies a domain $\Omega_h$ which can be described as a subgraph of a height profile function $h:\R\rightarrow [0,\infty)$, i.e., $\Omega_h:=\{\x:=(x,y)\in\R^2:\ 0<y<h(x)\}$. The energy functional is then given by 
\begin{eqnarray}
\label{eq:energy}
F(u,h):=\int_{\Omega_h} |\nabla u|^2\dxy+\int_{\R}\lt(\sqrt{1+h'^2}-1\rt)\dx,
\end{eqnarray}
where $u:\Omega_h\rightarrow\R$. For fixed volume $d>0$ of the film, we look for profile functions $h$ and associated displacement functions $u:\Om_h\rightarrow \R$ that minimize the total energy \eqref{eq:energy} 
subject to the constraints $\int_{\R}h\dx=d$ and $u(x,0)=e_0x$ for all $x\in\R$. The latter condition describes the crystallographic misfit between the substrate and the film, where  $e_0> 0$ measures its amplitude.  \\
The first term in \eqref{eq:energy} models the strain elastic energy in the film. Recall that we assume that there is a mismatch between the two crystal lattices, i.e., there is no stress-free configuration possible, and consequently, a strain is induced in the film while deposition.
 The second term in \eqref{eq:energy} models the extra surface energy due to the rearrangement of the atoms in  the film. All typical surface energy constants per unit length are normalized to one.
Let us notice that the functional $F$ bears many similarities with models for capillary surfaces \cite{mell:08,KimMel:12}. \\

We point out that, as noted before, in contrast to many previous works (see \cite{FFLM:07,FM:12,GolZwick}), we do not assume a periodic pattern of islands and do not restrict to a single island on a compact domain. The main difference is that in \eqref{eq:energy} the support of the height profile function $h$ may be unbounded,
which can lead to a loss of compactness for low energy sequences. A short comparison to the compact setting is given in Proposition \ref{rem:scalingcompls}. Many of our results, however, build on techniques developed in the works on compactly supported islands.\\

Let us make some comments on several simplifications built into the model. First, the displacement function $u\in H^1(\Omega_h)$ and the elastic energy term $|\nabla u|^2$ are scalar valued simplifications of a typical { geometrically} linear elastic energy density $W(U)=\mu|\frac{1}{2}(\nabla U+\nabla^T U)|^2$ 
for a displacement $U:\Omega_h\rightarrow\R^2$, $\mu$ being a typical elastic modulus. Based on the analysis in \cite{GolZwick}, we expect that the simplified energy contains, at least qualitatively, all relevant information. 
We note that the proofs of the scaling laws can be carried over to the elasticity setting, and are generalized to the three-dimensional setting
 in Section \ref{sec:3D}. Second, we assume that the domain occupied by the film can be described as a subgraph of the profile function $h$, which has the effect
 to prevent the formation of droplets or nanorings (see, e.g., \cite{ADFM:ADFM201302032}). Third, we do not take into account any plastic effects, such as misfit dislocations (see, e.g., \cite{NSM:91}). 
Finally, we consider only the stationary setting, and refer to \cite{FoFuLeMo:12,Pio:12,DalMFonLeo:12} for some recent results on the time evolution problem for the compact setting. \\

We consider two different types of approximations of the surface energy, namely for small  and large slopes $|h'|$. Many physical models are based on the assumption that for small volumes of the deposited film one expects small slopes of the film's profile function (see \cite{tersoff-tromp,SpeTer:10} or \cite{KimMel:12} where a similar simplification is used in the study of sliding liquid drops). This corresponds to the approximation (we ignore the factor $1/2$)
\begin{eqnarray}
\label{eq:smallslope}
\int_\R\lt(\sqrt{1+h'^2}-1\rt)\dx\sim \int_\R h'^2 \dx =:S_s(h).
\end{eqnarray}
If one expects, however, the formation of an island, the small slope approximation might not be appropriate anymore, and we compare it to the large slope approximation
\begin{eqnarray}
\label{eq:largeslope1}
\int_\R\lt(\sqrt{1+h'^2}-1\rt)\dx \sim \int_\R|h'|\dx=:S_\ell(h).
\end{eqnarray}
If we insert either of the approximations \eqref{eq:smallslope} or \eqref{eq:largeslope1} into \eqref{eq:energy}, then, due to the specific structure of the elastic energy term, we can rescale the problem to set $e_0=1$, i.e., we consider (see Section \ref{sec:smallslope} for a detailed derivation)
\begin{eqnarray}
\label{eq:energyresc}
 F_{s/\ell}(V):=\inf\left\{\int_{\Omh} |\nabla u|^2 \dxy + S_{s/\ell}(h)\ :\ h\in H^1(\R),\, { h \ge 0},\, \int_\R h \dx =V,\,  u(x,0)= x \right\}.
\end{eqnarray}
It turns out that in both cases there are two scaling regimes of the energy, 
namely 
(see Propositions \ref{prop:scaling} and \ref{prop:scalinglargeslop})
\[F_s(V)\sim \min\{V,\,V^{4/5}\}\mbox{\qquad and\qquad}F_\ell(V)\sim\min\{V,\,V^{2/3}\}. \]
 Heuristically, these scaling laws reflect the transition from a wetting regime in which the surface energy dominates and, consequently, the film forms
 a thin flat layer, to a regime in which a compactly supported island forms, { in which case the optimal energy comes from the competition}
between elastic and surface energy. We note that in this model, the surface energy prefers a flat layer, while the elastic energy favors oscillations. This is in contrast to many other physical situations where the surface energy typically favors compact shapes, and consequently minimizers exist for small volumes, but not for large volumes  (see, e.g., \cite{KM:2012,KM:2012-2,LuOtto:2012,GolNovRuf:13} for some recent works). The situation is opposite here: we prove that for large volumes, there always exists a minimizer, while for small volumes we prove non-existence of minimizers in the case of the small slope approximation (see Proposition \ref{prop:nonexist}). This is due to a loss of compactness of low-energy sequences, which corresponds to the wetting effect. 

To study more quantitatively the optimal shape of an island once it is formed, the limit $V\rightarrow\infty$ is considered.
If properly rescaled, the asymptotic shape turns out to be a parabola (in the case of the small slope approximation) and a rectangle (in the case of the large slope
 approximation).
{{We should stress the fact that even though it sheds some light on what can be (mathematically) predicted by our models, this asymptotic analysis does not tell much about the physics.
In fact, it corresponds to very large mismatch $e_0$ and/or large volume $d$, for which the model is not expected to be relevant anymore. Moreover, since in the original variables,
  minimizers tend to have large slopes, the small slope approximation is also questionable in this regime. We believe that adapting our analysis to the original functional \eqref{eq:energy} would lead to results similar to the one obtained in the large slope approximation.}}\\
 
The remaining part of the text is organized as follows. 
After setting the notation in Section \ref{sec:notation}, we first consider the small slope approximation. Some qualitative properties are derived in Section \ref{sec:firstprop}, and the scaling law for the minimal energy
is proven in Proposition \ref{prop:scaling}. The scaling law is then refined to show more quantitative results. More precisely, we show that there is a range of volumes $0<V<\overline{V}$, for which the minimal energy is exactly equal to $V$, and consequently, there does not exist a minimizer for $0 < V<\overline{V}$ (see Proposition \ref{prop:nonexist}). On the other hand, for volumes such that $F_s(V)<V$ (i.e., if $V > \overline{V}$), there always exists a smooth minimizer, which has compact and connected support and meets the substrate at zero angle (see Proposition \ref{existminss} and Theorem \ref{thm:reg}). Regularity properties and estimates on the support and its maximal height are provided (see Sections \ref{sec:reg} and \ref{sec:morequal}). { The regularity of a minimizer (see Section \ref{sec:reg}) is shown using 
arguments from \cite{FFLM:07} with minor changes required due to a different form of the surface energy. Though the proof follows a standard approach, for the sake of completeness we include it in the paper.} Finally, the asymptotic behavior for large $V\rightarrow\infty$ is studied. It is shown, that, when properly rescaled, minimizers converge to a parabola, and away from a boundary layer, this convergence occurs at an  exponential rate in the $L^2$-topology.

Subsequently, in Section \ref{sec:largeslope}, the large slope approximation \eqref{eq:largeslope1} is considered. Note that this approximation comes along with a loss of regularity of admissible profile functions $h$, and we consider the relaxation following \cite{FFLM:07}. We prove the scaling law of the minimal energy $F_\ell$ (see Proposition \ref{prop:scalinglargeslop}), and show that in the regime $F_\ell(V)<V$ there always exists a minimizer with connected support (see Proposition \ref{prop:existminls}). If properly rescaled, a sequence of minimizing profiles converges, away from a boundary layer, to a rectangular shape for large volumes at an exponential rate in the $L^1$-topology (see Proposition \ref{prop:exp1}). Recall that for the small slope approximation, non-existence of minimizers at small volumes is due to the fact that $F_s(V)=V$. For the large slope approximation, we only get the weaker result $F_\ell(V)/V\rightarrow 1$ as $V\rightarrow 0$ (see Proposition \ref{prop:smallvollargeslope}). 

Finally, in Section \ref{sec:3D}, the three-dimensional setting is considered, and the scaling laws for both types of approximations are discussed.

\section{Notation and preliminary results}\label{sec:notation}

In this section, we set the notation and collect some results that will be used later. Throughout the text we denote by $C$ and $c$ constants that may vary from expression to expression. The symbols $\sim$, $\ges$, $\les$ indicate estimates that hold up to a constant. For instance, $f\les g$ denotes the existence of a constant $C>0$ such that $f\le Cg$.
For $\Om\subset \R^2$, we denote by $\mathcal{H}^1\left(\Om\right)$ its one-dimensional Hausdorff measure, and by $\left|\Om\right|$ its two-dimensional Lebesgue measure. When it exists, we will denote by $\nu$ its inward normal. Given two sets $A$, $B\subset \R^2$, we define their Hausdorff distance as $d_{\mathcal{H}}(A,B):=\inf\{r>0 \, :\, A\subset N(B,r) \textrm{ and } B\subset N(A,r)\}$, where $N(A,r):=\{x\in \R^2 \ :\ d(x,A)<r\}$, and $d(x,A)$ denotes the distance from $x$ to $A$.    
 Given a vector $\x:=\left(x,y\right)\in \R^2$, we denote by
 $\left|\x\right|:=\left(x^2+y^2\right)^{1/2}$ its Euclidean norm. 

We will use the following rescaling property for functions on rectangles.
\begin{lemma1}\label{lem:resc}
If $u\in H^1\left(\left[0,\ell\right]\times\left[0,L\right];\R\right)$, with $u\left(x,0\right)=x$, then letting  $v\left(x,y\right):= \frac{1}{\ell} u\left(\ell x, \ell y\right)$, there holds
\begin{equation}\label{estimenergyfonda}
\int_{\left[0,\ell\right]\times\left[0,L\right]}|\nabla u|^2 \dxy   = \ell^2 
\int_{\left[0,1\right]\times\left[0,L/\ell\right]}| \nabla v|^2 \dxy. \end{equation}
As a consequence, for every $\lambda>0$, there exists $C(\lambda)>0$ such that
 \begin{equation}\label{estimenergyfondabis}
\int_{\left[0,\ell\right]\times\left[0,\lambda \ell\right]}|\nabla u|^2 \dxy   \ge C(\lambda) \ell^2. \end{equation}
\end{lemma1}

The following lemma describes the behavior of the elastic energy for small thickness of the film. It can be seen as a simple special case of a dimension reduction argument (see \cite{LDR:95,AnBaPer:94}). 
{{A proof for the more complicated case of vector-valued functions}}
can be found in \cite{GolZwick}.
 
\begin{lemma1}\label{lemLeDreRaou}
 There holds
 \begin{equation}\label{LeDretRaoult}\lim_{\eps\to 0^+} \min_{u(x,0)=x} \frac{1}{\eps} \int_{\left[0,1\right]\times\left[0,\eps\right]} \left|\nabla u\right|^2\dxy = 1.\end{equation}
\end{lemma1}
\begin{remark1}
 The analogous statements hold for typical elastic energy functionals for deformations $U:\R^2\rightarrow\R^2$, i.e.,
 \begin{equation}\label{LeDretRaoult2}\lim_{\eps\to 0^+} \min_{U(x,0)=(x,0)} \frac{1}{\eps} \int_{\left[0,1\right]\times\left[0,\eps\right]} W(\nabla U)\dxy = 1\end{equation}
if $W(\nabla U)=|\nabla U|^2$ or $W(\nabla U)=|\frac{1}{2}(\nabla U+\nabla^T U)|^2$ (see \cite{GolZwick}). This allows to carry over the qualitative results to the linear elasticity setting.
\end{remark1}
\section{Small slope approximation}\label{sec:smallslope}
In this section we consider the small slope approximation $\sqrt{1+h'^2}-1 \sim h'^2$ and study 
\begin{eqnarray}\label{funcssor}
 F_s(e_0,d)&:=&\inf\left\{\int_{\Omh} |\nabla u|^2\, \dxy \, +\, \int_{\R} h'^2 \dx\ :\ {{h\in H^1(\R),\ h\geq 0, \ u\in H^1(\Omega_h),\ }}\right. \nonumber\\
&&\quad \left.\int_\R h\, \dx=d,\,  u(x,0)=e_0 x {{\ \text{if\ }x\in\supp h}}\right\}.
\end{eqnarray}
{{One main difference between the model considered here and the related models on compact domains (see \cite{FFLM:07,FM:12,GolZwick}) is that
it behaves well under rescaling as shown by the following lemma.}}
\begin{lemma1}\label{rescale}
 For $h\in {{H}}^1(\R)$, $ h \ge 0$, $u\in H^1(\Omh)$, and $\lambda>0$, letting $u_\lambda(x,y)=\frac{1}{\lambda} u(\lambda x, \lambda y)$ and $h_\lambda(x)=\frac{1}{\lambda} h(\lambda x)$, we have
\[  \int_{\Om_{h}} |\nabla u|^2\, \dxy\, +\, \int_{\R} h'^2\, \dx= \lambda^2\int_{\Om_{h_\lambda}} |\nabla u_\lambda|^2\,\dxy \,  +\, \lambda \int_{\R} h_\lambda'^2 \, \dx\]
and $\int_{\R} h\, \dx= \lambda^2 \int_{\R} h_\lambda \, \dx$.
\end{lemma1}
{{By rescaling, we can eliminate one of the two parameters $d$ or $e_0$. We will renormalize such that $e_0=1$:
\begin{prop1}
 Let $e_0>0$ and $d>0$. Set $V:=e_0^4d$ and
\begin{eqnarray}
 \label{funcss}
 F_s(V):=\inf\left\{\int_{\Omh} |\nabla u|^2\,\dxy +\, \int_{\R} h'^2\, \dx\ :\ (u,h)\in\mathcal{A}_V \right\},
\end{eqnarray}
where the set of admissible pairs is given by
\begin{eqnarray}
 \label{eq:admset}
\mathcal{A}_V:=\left\{ (u,h):\ h\in H^1(\R),\ h\geq 0, \ u\in H^1(\Omega_h),\ \int_\R h\, \dx=V,\,  u(x,0)= x \ \text{if\ }x\in\supp h\right\}.
\end{eqnarray}
Then 
\[ F_s(e_0,d)=\frac{1}{e_0^2} F_s(V).\]
\end{prop1}
\begin{proof}
 Let $(h,u)$ be admissible for $F_s(e_0,d)$, and let $\lambda:=\frac{1}{e_0^2}$. Using the notation of Lemma \ref{rescale}, set
\begin{eqnarray}\label{eq:e0dV}
 \tilde{h}(x):=h_\lambda(x),\qquad \text{and\quad}\tilde{u}(x,y):=\frac{1}{e_0}u_\lambda(x,y).
\end{eqnarray}

Then $(\tilde{u},\tilde{h})\in\mathcal{A}_V$ and, by Lemma \ref{rescale},
\[\int_{\Omega_h}|\nabla u|^2\dxy+\int_\R h'^2\dx=\frac{1}{e_0^2}\left(\int_{\Omega_{\tilde{h}}}|\nabla \tilde{u}|^2\dxy+\int_{\R}\tilde{h}'^2\dx'    \right).\]
Since \eqref{eq:e0dV} induces a bijective correspondence between the admissible pairs for $F_s(e_0,d)$ and $\mathcal{A}_V$, this proves the assertion.
\end{proof}

}}


In this section we study the problem \eqref{funcss}. We first prove the scaling law of the optimal energy. It is shown that there exists a critical volume $\overline{V}>0$ such that  for volumes $0 < V < \overline{V}$, we have $F_s(V)=V$, which leads to the non-existence of minimizers. We also prove that for $V > \overline{V}$, we have $F_s(V)<V$ and there exists a compact connected smooth minimizer of \eqref{funcss}, which has zero contact angle with the substrate. Moreover, we provide estimates on the size of the support of this island together with 
estimates on its maximal height. We finally investigate the large volume limit and prove that, when suitably rescaled,  the minimizers converge to a parabola and that away from a boundary layer this convergence is exponentially fast in the (strong) $L^2$-topology. {{We point out that we investigate the asymptotic behavior mostly for mathematical reasons to better understand the shapes of minimizers. The physical model, and in particular the small slope approximation implemented here, is expected to give more reliable results in the case of small volumes.}}

\subsection{First properties of the minimization problem}\label{sec:firstprop}

For every $V>0$ and every $(u,h) \in\mathcal{A}_V$, we set
\[E_V(u,h):=\int_{\Omh}|\nabla u|^2 \dxy \qquad \textrm{and} \qquad S_V(h):= \int_\R h'^2 \dx.\]
When it is clear from the context, we will often drop the explicit dependence on $(u,h)$. 
As a simple consequence of Lemma \ref{rescale}, we have the following important property of $F_s$. 
\begin{prop1}\label{concavss}
 {{$F_s:(0,\infty)\rightarrow\R$}} is 
concave 
and 
thus locally Lipschitz continuous.
\end{prop1}
\begin{proof}
 Let $V>0$ be fixed. Then, for every $V_0 > 0$ and for every admissible competitor $(u,h)    { \in\mathcal{A}_{V_0}}$, by Lemma \ref{rescale}, $F_s(V) \le \frac{V}{V_0} E_{V_0}(u,h) +\left(\frac{V}{V_0}\right)^{1/2} S_{V_0}(h)$. Hence
\begin{equation}\label{concavF}
F_s(V)=\inf_{V_0>0} \inf_{ (u,h) \in \mathcal{A}_{V_0}}\, \left( \frac{V}{V_0} E_{V_0}(u,h) +\left(\frac{V}{V_0}\right)^{1/2}S_{V_0}(h)\right),
\end{equation}
and since for every $V_0 > 0$ and $(u, h) {{\in\mathcal{A}_{V_0}}}$, the function $V\mapsto \frac{V}{V_0} E_V(u,h) +\left(\frac{V}{V_0}\right)^{1/2}S_V(h)$ is concave, $F_s$ is the infimum of concave functions and therefore also concave.
\end{proof}
The scaling behavior from Lemma \ref{rescale} is typical for various discrete and continuous models for epitaxial growth, and thus, similar properties hold for a large class of models (see also \cite{FonPraZwi:14,SpeTer:10}). Using the same rescaling we obtain the following result.
\begin{prop1}\label{connect}
 For every $V>0$, {{ if 
\[F_s(V)=\min\left\{\int_{\Omh} |\nabla u|^2 \dxy +\, \int_{\R} h'^2 \dx : (u,h)\in\mathcal{A}_V
\right\} =E_V(u^\ast,h^\ast)+S_V(h^\ast),\]
i.e., if the minimum is attained, then $\{h^\ast>0\}$ is connected.}}
\end{prop1}
{{
\begin{proof}
Fix $V>0$, and suppose that there exists a minimizer. It suffices to show that for every $(u,h)\in\mathcal{A}_V$, for which $\{h>0\}$ is not connected, there exists
 $(\overline{h},\overline{u})\in\mathcal{A}_V$ with lower total energy. Assume that $h=h^{(1)}+h^{(2)}$ with $h^{(1)}\geq 0$, $h^{(2)}\geq 0$,
 $\{h^{(1)}>0\}\cap\{h^{(2)}>0\}=\emptyset$, $0<\int_{\R}h^{(1)}\dx=V_1<V$, $\int_{\R}h^{(2)}\dx=V-V_1$, and $h^{(1)},h^{(2)} \in H^1(\R)$. Up to translation, we can further assume that $h^{(1)}$ and $h^{(2)}$ were chosen such that 
\begin{equation}\label{separcomp}
 \{h^{(1)}>0\} \subset \R^- \qquad \textrm{and} \qquad \{h^{(2)}>0\} \subset \R^+.
\end{equation}
   We set $u^{(1)}:=u|_{\Omega_{h^{(1)}}}$ and $u^{(2)}:=u|_{\Omega_{h^{(2)}}}$. Then $(u^{(i)},h^{(i)})\in\mathcal{A}_{V_i}$ for $i=1,2$. We note that
\[E_V(u,h)+S_V(h)= E_{V_1}(u^{(1)},h^{(1)})+S_{V_1}(h^{(1)})+ E_{V_2}(u^{(2)},h^{(2)})+S_{V_2}(h^{(2)}). \]
We now build a competitor as follows: Consider the two components separately and rescale them. Precisely, for $0\leq \mu\leq V$, set $\lambda_1:=\sqrt{\frac{V_1}{\mu}}$ and $\lambda_2:=\sqrt{\frac{V_2}{V-\mu}}$, and,
 using the notation of Lemma \ref{rescale} consider 
\[h_\mu:=h^{(1)}_{\lambda_1}+h^{(2)}_{\lambda_2}(\cdot-\tau_\mu), \]
where $\tau_\mu$ is such that $\{h^{(1)}_{\lambda_1}>0\}\cap\{h^{(2)}_{\lambda_2}(\cdot-\tau_\mu)>0\}=\emptyset$ (which exists thanks to \eqref{separcomp}), and the associated $u_\mu$. Note that $h=h_{V_1}$, and, by Lemma \ref{rescale}, $\int_{\R}h_{\mu}\dx=\frac{\mu}{V_1}V_1+\frac{V-\mu}{V_2}V_2=V$ for every $\mu\in[0,V]$. Hence $(h_\mu,u_\mu)\in\mathcal{A}_V$ for every $\mu\in [0,V]$. 
Let \[f(\mu):= E_V(h_\mu,u_\mu)+S_V(h_\mu),\]
 so that  by Lemma \ref{rescale},
\begin{eqnarray*}
 f(\mu)=\frac{\mu}{V_1}E_{V_1}(h^{(1)},u^{(1)})+\sqrt{\frac{\mu}{V_1}}S_{V_1}(h^{(1)},u^{(1)})+\frac{V-\mu}{V_2}E_{V_2}(h^{(2)},u^{(2)})+\sqrt{\frac{V-\mu}{V_2}}S_{V_2}(h^{(2)},u^{(2)}).
\end{eqnarray*}
Since $S_{V_1}(h^{(1)},u^{(1)})>0$ and $S_{V_2}(h^{(2)},u^{(2)})>0$, the function $f$ is strictly concave. Therefore, it attains its minimum at the boundary, that is, at $\mu=0$ or at $\mu=V$. This shows that there is a configuration with strictly lower energy than $(u,h)$ and proves that the minimizer must be connected. 
\end{proof}

}}
 Using a different rescaling we obtain the following:
\begin{lemma1}\label{internvar}
If $V>0$, then 
\[F_s(V)=\inf\left\{\int_{\Omh} |\nabla u|^2\dxy + \int_{\R} h'^2 \dx : 
{{(u,h)\in\mathcal{A}_V}};\  {{\partial_yu\equiv0\text{\ or\ }}}
 \int_{\Omega_h} (\partial_y u)^2 \dxy = \frac{3}{4} \int_\R h'^2 \dx
\right\}.\]
\end{lemma1}

\begin{proof}
{{For $(u,h)\in\mathcal{A}_V$ consider the equivalence class $\{(u_\lambda,h_\lambda):\ \lambda>0\}\subset\mathcal{A}_V$ given by}}
the anisotropic rescaling $u_\lambda(x,y)=\frac{1}{\lambda} u(\lambda x, \frac{1}{\lambda} y)$ and $h_\lambda(x)=\lambda h(\lambda x)$. 
{{Then}}
\[\int_{\Om_{h_\lambda}} |\nabla u_\lambda|^2\dxy + \int_{\R} h_\lambda'^2 \dx=\int_{\Omh} \left((\partial_x u)^2+\frac{1}{\lambda^4} (\partial_y u)^2\right) \dxy + \lambda^3\int_{\R} h'^2 \dx.\]
{{Suppose that $\partial_yu\not\equiv 0$. Since for every $(u,h)\in\mathcal{A}_V$ we have $\int_\R
h'^2\dx>0$, within one equivalence class, the energy is minimized for $\lambda=\frac{4}{3}\frac{\int_{\Omega_h} (\partial_y u)^2 \dxy}{\int_\R h'^2 \dx}>0$. Therefore, from each such equivalence class, only the element with $\int_{\Omega_h} (\partial_y u)^2 \dxy = \frac{3}{4} \int_\R h'^2 \dx$ is relevant for the infimum.
}}
\end{proof}
Following \cite{BC:02,FFLM:07,GolZwick}, we prove the lower semicontinuity of the energy and density of Lipschitz configurations.
\begin{prop1}\label{densityss}
 For every sequence $(h_n,u_n){{\in\mathcal{A}_V}}$ 
{{with}} $\sup_n E_{V}(u_n,h_n)+S_V(h_n) \le C$, there exists $(u,h)$ {{with $h\in H^1(\R), h \ge 0$ and $u\in H^1(\Omega_h)$}} such that up to a subsequence, $h_n$ converges to $h$ in $L^\infty_{\text{loc}}(\R)$,  $u_n$ converges weakly in 
$H^1_{loc}(\Omega_h)$ to $u$, $u(x,0)=x$, and
\[\varliminf_{n\to +\infty} E_{V}(u_n,h_n)+S_V(h_n)  \ge \int_{\Omega_h} |\nabla u|^2 \dxy+\int_{\R} h'^2 \dx.\]
Moreover, if $\{h_n\}$ is tight, then $h_n\rightarrow h$ in $L^1(\R)$, and $\int_{\R}h\dx=V$. 

Conversely, for every $(u,h){{\in\mathcal{A}_V}}$ with 
$E_{V}(u,h)+S_V(h)<\infty$, there exists a sequence $(h_n,u_n){{\in\mathcal{A}_V}}$ 
{{such that}} $h_n$ {{is}} Lipschitz continuous with bounded support and such that  $h_n$ converges to $h$ in $L^\infty(\R)$,  $u_n$ converges weakly in 
$H^1_{loc}(\Omega_h)$ {{to $u$,}} and
\[\varlimsup_{n\to +\infty} E_{V}(u_n,h_n)+S_V(h_n)  \le E_{V}(u,h)+S_V(h).\]
\end{prop1}
\begin{proof}
 Let $(h_n,u_n){{\in\mathcal{A}_V}}$ be a sequence such that $\sup_n E_{V}(u_n,h_n)+S_V(h_n) \le C$.
From the compact embedding of $H^1_{loc}(\R)$ in $L^{\infty}_{loc}(\R)$, we get that (up to a subsequence),
 $h_n$ converges in $L^{\infty}_{loc}(\R)$ to some (continuous) function $h$. This implies the local Hausdorff convergence of $\R^2\backslash \Omega_{h_n}$ to $\R^2\backslash \Omega_h$, from which as in \cite{BC:02,FFLM:07} we infer the existence of a function $u{{\in H^1(\Omega_h)}}$ with $u(x,0)=x$ such that 
$u_n$ converges weakly in $H^1_{loc}(\Omega_h)$, and as in \cite{BC:02,FFLM:07} we obtain the lower semicontinuity
\[\varliminf_{n\to +\infty} E_{V}(u_n,h_n)+S_V(h_n)  \ge\int_{\Omega_h} |\nabla u|^2 \dxy+\int_{\R} h'^2 \dx.\]
If $\{h_n\}$ is tight, then $h_n\rightarrow h$ in $L^1(\R)${{, which implies that $\int_{\R}h\dx=V$}}. 

Let now $(u,h)  \in \mathcal{A}_V$ be such that 
$E_{V}(u,h)+S_V(h)<\infty$. It is readily seen that we can approximate  $h$ 
from below 
by compactly supported height profiles so that we will assume from now on that $h$ itself is compactly supported.
It {{suffices}} to approximate $h$ by a sequence of Lipschitz functions $h_n$ with $0\le h_n\le h$ and $S_V(h_n)\le S_V(h)$. We refer the reader to  \cite{BC:02,FFLM:07,GolZwick} for the treatment of the volume constraint. Following \cite{BC:02}, for $n\in \N$ and $x\in \R$, we define 
\[h_n(x):=\inf_{x' \in \R} \left( h(x')+n|x-x'| \right)\]
to be the Yosida transform of $h$. {{The latter}} is an $n-$Lipschitz function which satisfies $0\le h_n\le h$. As proven in \cite{BC:02}, $\Omega_{h_n}$ converges to $\Omega_h$ in the Hausdorff topology. Since $h_n$ and $h$ are continuous functions, the set $\{h_n< h\}$ is open and thus made of a countable union of disjoint intervals $(a_k,b_k)$, $k\in \N$. On each of these intervals (see \cite{BC:02}),
   
\[h_n(x)=\min\{h(a_k)+ n|x-a_k|, h(b_k)+ n|x-b_k|\}\]
so that $\int_{a_k}^{b_k} h_n'^2 \dx \le \int_{a_k}^{b_k} h'^2 \dx$. From $h'=h'_n$ a.e. on $\{h_n=h\}$, we obtain that $S_V(h_n)\le S_V(h)$ and that $h_n$ converges in $L^\infty$ to $h$.
\end{proof}

We now prove an interpolation inequality which will be useful later.
\begin{prop1}\label{prop:interpol}
 For every $h{{\in H^1(\R)}}$, {{we have}}
\begin{equation}\label{interpolationineq}
 \|h\|_{L^\infty(\R)} \le \left( \frac{9}{16} \right)^{1/3} \left( \int_\R |h(x)| \dx \right)^{1/3} \left( \int_\R h'(x)^2 \dx \right)^{1/3}.
\end{equation}

\end{prop1}

\begin{proof}
{{Without loss of generality, we may assume that $h\geq 0$ (otherwise consider $|h|$). By rescaling the dependent and the independent variables, we may assume that $\|h\|_{L^\infty}=1$ and $\int_\R h \dx= 1$. 
Indeed, suppose that the inequality holds for some $h\in H^1$. For $M>0$ and $\lambda>0$ consider the rescaled function $\tilde{h}(x)=Mh(\lambda x)$. Then, by changing
 variables, $\tilde{x}=\frac{x}{\lambda}$, we get that \eqref{interpolationineq} holds for $\tilde{h}$ since
\begin{eqnarray*}
 \|\tilde{h}\|_{L^\infty}&=&M\|h\|_{L^\infty} \leq M\left( \frac{9}{16} \right)^{1/3} \left( \int_\R |h(x)| \dx \right)^{1/3} \left( \int_\R h'(x)^2 \dx \right)^{1/3}\\
&=& M\left( \frac{9}{16} \right)^{1/3} \left( \frac{1}{M}\int_\R |(Mh)(x)| \dx \right)^{1/3} \left(\frac{1}{M^2} \int_\R (Mh)'(x)^2 \dx \right)^{1/3}\\
&=& \left( \frac{9}{16} \right)^{1/3} \left(\lambda \int_\R |(Mh)(\lambda \tilde{x})| \, \mathrm{d} \tilde{x} \right)^{1/3} \left( \lambda\int_\R \frac{1}{\lambda^2}(Mh)'(\lambda\tilde{x})^2 \, \mathrm{d} \tilde{x} \right)^{1/3}\\
&=& \left( \frac{9}{16} \right)^{1/3} \left( \int_\R |\tilde{h}(x)| \dx \right)^{1/3} \left( \int_\R \tilde{h}'(x)^2 \dx \right)^{1/3}.
\end{eqnarray*}
By translation and symmetric decreasing rearrangement (see \cite{liebloss}), we can further restrict ourselves to functions $h$ with $h(0)=\sup h=1$, which are even and non-increasing on $[0,+\infty)$. Let
 }}


\begin{equation*}
 g(x) := \begin{cases} \left( 1 - \frac{|x|}{x_0}\right)^2, & \text{if\ }x \in (-x_0,x_0), \\ 0, & \textrm{else}, \end{cases}
\end{equation*}
where $x_0>0$ is chosen so that the volume constraint is satisfied, i.e., $x_0=3/2$. {{Let us}} prove that $g$ is the minimizer {{of $\int_{\R}h'^2\dx$ in the class
\[\mathcal{M}:=\left\{h\in H^1(\R):\ h\geq 0, \ \|h\|_{L^\infty}=h(0)=1, \ \int_{\R}h\dx=1,\  h(x)=h(-x),\ h\text{\ non-increasing on }[0,\infty)\right\}. \] Let $h\in\mathcal{M}$. 
Then
\begin{align*}
 \int_0^{+\infty} h'^2 \dx &=\int_0^{+\infty} g'^2\dx+ 2\int_0^{+\infty} g'(h-g)' \dx +\int_0^{+\infty} |(h-g)'|^2 \dx\\
&=\int_0^{+\infty} g'^2\dx- 2\int_0^{x_0} g''(h-g) \dx +\int_0^{+\infty} |(h-g)'|^2 \dx
\end{align*}
where we used integration by parts and the fact that $g'(x_0)=0$, $h(0)=g(0)=1$. Since on $[0,x_0]$, $g''=\frac{2}{x_0^2}$ and since $\int_0^{+\infty} (h-g) \dx=0$, we further obtain that
\[\int_0^{x_0} g''(h-g) \dx=\frac{2}{x_0^2}\int_0^{x_0} (h-g) \dx=-\frac{2}{x_0^2}\int^{+\infty}_{x_0} h \dx\]
so that 
\[\int_0^{+\infty} h'^2 \dx =\int_0^{+\infty} g'^2\dx+ \frac{4}{x_0^2}\int^{+\infty}_{x_0} h \dx +\int_0^{+\infty} |(h-g)'|^2 \dx\ge \int_0^{+\infty} g'^2\dx.
\]
By symmetry, this shows that $g$ minimizes $\int_{\R}h'^2\dx$ in $\mathcal{M}$. 
Using that $\int_0^{x_0} g'^2 \dx=8/9$, we obtain~\eqref{interpolationineq}.
}}

\end{proof}
\begin{remark1}\label{rem:interpol}
 In terms of the energy, \eqref{interpolationineq} can be rephrased as
\[\|h\|_{L^\infty(\R)} \le  \left(\frac{9}{16} \right)^{1/3} V^{1/3}  S_V(h)^{1/3}.\]
\end{remark1}

\subsection{Scaling law}

In this section, we prove the following scaling law for the energy.
\begin{prop1}\label{prop:scaling}
 There exists a 
constant $\c{{>0}}$ such that for every $V>0$,
\[\c \min\{V,\, V^{4/5}\}\le F_s(V)\le  \c^{-1} \min\{V,V^{4/5}\}.\]
Moreover, {{there exists $C>0$ with the following property: If}} $V$ is large enough, and $(u,h){{\in\mathcal{A}_V}}$ {{with}} $E_V(u,h)+S_V(h) \le \frac{1}{\c} V^{4/5}$, then $\max h \ge C V^{3/5}$. 
\end{prop1}
\begin{proof}
{{We prove the upper bound first. For that, we have to construct two elements from $\mathcal{A}_V$ for $V>0$. First, for $N\in\N$, set 
\[h_N(x):=\begin{cases}
           \frac{V}{N^2}x+\frac{V}{N},&\text{if\ } -N\leq x\leq 0,\\
-\frac{V}{N^2}x+\frac{V}{N},&\text{if\ }0\leq x\leq N,\\
0,&\text{if\ }|x|\geq N,
          \end{cases}
 \]
and set $u_N(x,y):=x$ if $(x,y)\in\Omega_{h_N}$. Then $(u_N,h_N)\in\mathcal{A}_V$, and $E_V(h_N,u_N)+S_V(h_N)=V+\frac{2V^2}{N^3}\leq 2V$ for $N$ large enough. Note that we even have
\begin{eqnarray}\label{eq:Vupper}
 E_V(h_N,u_N)+S_V(h_N)\rightarrow V\text{\quad as\ }N\rightarrow\infty.
\end{eqnarray}
For the other regime, set
\[\tilde{h}(x):=\begin{cases}
                 V^{1/5}x+V^{3/5},&\text{if\ }-V^{2/5}\leq x\leq 0,\\
-V^{1/5}x+V^{3/5},&\text{if\ }0\leq x\leq V^{2/5},\\
0,&\text{if\ }|x|\geq V^{2/5},
                \end{cases}
 \]
and let $\tilde{u}$ }}
be the restriction to $\Omega_h$ of 
\[{{\tilde{u}}}(x,y):=\begin{cases}
        x(1-V^{-2/5}y),& \mbox{\quad if\ } y\le V^{2/5},\\
0,&\mbox{\quad if\ } y\ge V^{2/5}.  
         \end{cases}\]
{{
Then $(\tilde{u},\tilde{h})\in\mathcal{A}_V$, and 
\begin{eqnarray*}
 E_V(\tilde{u},\tilde{h})+S_V(\tilde{h})&\leq&\int_{-V^{2/5}}^{V^{2/5}}\int_0^{V^{2/5}}(1-V^{-2/5}y)^2+V^{-4/5}x^2\dxy+2\int_0^{V^{2/5}}V^{2/5}\dx\\
&=&\lt(\frac{10}{3}+2\rt)V^{4/5}.
\end{eqnarray*}
}}
\noindent 
For the lower bound, thanks to Proposition \ref{densityss}, we can assume that $h$ is a Lipschitz function with $\int_\R h \dx=V$ and $\{h>0\}$ compact and connected {{(otherwise consider each of the (possibly infinitely many) connected components separately)}}. Let then $u\in H^1(\Omh)$ be the minimizer of the Dirichlet energy {{in $\Omega_h$}} with $u(x,0)=x$. Let $x_1\in\R$ be such that $h(x_1)>0$, and let $t_1 > 0$ be the maximal $t>0$ such that the square $[x_1,x_1+t]\times[0,t]$ is below the graph of $h$, i.e., $[x_1,x_1+t]\times[0,t] \subset \overline{\Omh}$. Observe that the maximality of $t_1$ implies the existence of a point $\bar x_1 \in [x_1,x_1+t_1]$ such that $h(\bar x_1) = t_1$. Now let
\begin{equation*}
 V_1 := \int_{x_1}^{x_1+t_1} h \dx\, , \quad E_1:=\int_{x_1}^{x_1+t_1} \int_0^{h(x)} |\nabla u|^2 \dxy, \quad \textrm{ and } \quad S_1:=\int_{x_1}^{x_1+t_1} h'^2 \dx.
\end{equation*}
We want to show that $E_1 + S_1 \gtrsim \min \{V_1,V_1^{4/5}\}$. By Lemma~\ref{lem:resc} ,  $E_1 \ge \int_{[x_1,x_1+t_1]\times [0,t_1]} |\nabla u|^2 \dx \ge C t_1^2$, hence $E_1 \ge \frac{C}{2} V_1$ provided that $V_1 \le 2t_1^2$. 

Let us now assume that $V_1 > 2t_1^2$. Since $\max_{[x_1,x_1+t_1]} h \ge \frac{V_1}{t_1}$, we have 
\[S_1 \ge \frac{1}{t_1} \left(\max_{[x_1,x_1+t_1]} h -\min_{[x_1,x_1+t_1]} h \right)^2 \ge 
\frac{1}{t_1} \left(\frac{V_1}{t_1} - t_1\right)^2 = \frac{1}{t_1} \left(\frac{V_1}{2t_1} + \frac{V_1}{2t_1} - t_1\right)^2 \ge \frac{V_1^2}{4t_1^3}.\] Since $E_1 \ge C t_1^2$, we obtain
\begin{equation}\nonumber
 E_1 + S_1 \ge C t^2_1 + \frac{V_1^2}{4t_1^3} \ge C V_1^{4/5},
\end{equation}
where the last inequality follows from  Young's inequality. 

We have thus shown that the energy in $[x_1,x_1+t_1]$ is bounded from below by $\c \min\{V_1,V_1^{4/5}\}$. We define iteratively  $x_{2} := x_1+t_1$, $x_{i+1}:=x_1+\sum_{k=1}^it_k$, and repeat the process in each interval $[x_i,x_i+t_i]$, 
and similarly in the opposite direction (i.e., going to the left) starting at $x_1$.
Since $h$ is Lipschitz we cover with this procedure the whole set $\{h>0\}$, and the lower bound follows.

It remains to show that 
{{if $V$ is large enough, and $(u,h)\in\mathcal{A}_V$ with $E_V(u,h)+S_V(h) \le \frac{1}{\c} V^{4/5}$, then $\max h \ge C V^{3/5}$ for some $C>0$ independent of $V$.}} 
Having $t_i$ and $V_i$ constructed in the previous part of the proof, let us assume that $V_1$ is the largest among all $\{V_i\}$. Then we have
\begin{multline*}
 {{\frac{V^{4/5}}{\c}}} \ge E_V(u,h) + S_V(h) \ge \sum_{{{i}}} \c\min\{V_i,V_i^{4/5}\}\\
 \ge \c \min\{1,V_1^{-1/5}\} \sum_{{{i}}} V_i = \c \min \{1,V_1^{-1/5}\} V,
\end{multline*}
which gives $c_0^{10} V\le \max\{1,V_1\}$. Assume now that $V \ge \c^{-10}$ so that $V_1\ge c_0^{10} V$. If we set $M:=\max_{[x_1,x_1+t_1]} h$, then as before $M \ge V_1/t_1$, which implies $t_1 \ge V_1/M$. Since $E_1 \ge C t_1^2$, we have
\[
 V^{4/5}/\c \ge E_1 \gtrsim t_1^2 \ge (V_1/M)^2 \ge \c^{20} V^2 M^{-2},
\]
which implies $\max h\geq M \gtrsim V^{3/5}$. 
\end{proof}

\begin{remark1}
%
%
%
%
%
%
Using the scaling law (see Proposition \ref{prop:scaling}) and \eqref{interpolationineq}, we find that $\max h\les V^{3/5}$, and so the previous proposition implies that $\max h \sim V^{3/5}$ 
for sufficiently large $V$.
Using this we see that the size of the support of $h$ is at least of order $\max\{1,V^{2/5}\}$. 
\end{remark1}
In terms of the original parameters $e_0$ and $d$, the scaling law reads as follows.
\begin{prop1}
 There exists a positive constant $c$ such that for every $V>0$,
\[c \min\{e_0^2d,\, d^{4/5} e_0^{6/5}\}\le F(V)\le \frac{1}{c} \min\{e_0^2d,\, d^{4/5} e_0^{6/5}\}.\]
\end{prop1}
\begin{remark1}
 In the original coordinates, for $e_0^4 d \ges 1$, the typical island is of height $e_0^{2/5}d^{3/5}$ and of width  $\left(\frac{d}{e_0}\right)^{2/5}$.
\end{remark1}

\subsection{Existence and non-existence of minimizers}

Let us start by studying the non-existence {{case for small $V$}},
i.e., the wetting regime. For this, we 
prove that for a non-trivial range of volumes {{$0<V\leq\overline{V}$, }}the infimum of the energy $F_s(V)$ is equal to $V$. 
\begin{prop1}\label{prop:nonexist}
 There exists $\overline{V}>0$ such that for every $0 < V\le \overline{V}$ we have $F_s(V)=V$. As a consequence, there exists no minimizer of \eqref{funcss} for $0 < V < \overline{V}$.
\end{prop1}
\begin{proof}
{{Recall that by the construction of flat layers in the proof of the scaling law (see \eqref{eq:Vupper}), we have $F_s(V)\leq V$. It remains to show the reverse inequality, i.e., that for every $(u,h)\in\mathcal{A}_V$, we have $E_V(u,h)+S_V(h)\geq V$ if $V$ is small enough. By density, it suffices to consider $(u,h)\in\mathcal{A}_V$ for which $h$ is Lipschitz continuous (see Proposition \ref{densityss}), and $u$ being the minimizer of the Dirichlet energy in $\Omega_h$ subject to the boundary condition $u(x,0)=x$ on $\{h>0\}$. In particular, $-\Delta u=0$ in $\Omega_h$. 
}}
Testing the Laplace equation with $u(x,y) - x$, we obtain
\begin{equation*}
 0 = \int_{\Omega_h} \div (\nabla u) \cdot (u-x) \dxy = -\int_{\Omega_h} \lt(|\nabla u|^2 - \partial_x u\rt) \dxy + \int_{\partial \Omega_h} \partial_\nu u \cdot (u-x) \ud \mathcal{H}^1.
\end{equation*}
The boundary integral vanishes, because $\partial_\nu u = 0$ if $y = h(x)$, and $u(x,0)=x$. Hence 
\begin{equation*}
 \int_{\Omega_h} |\nabla u|^2 \dxy = \int_{\Omega_h} \partial_x u \dxy,
\end{equation*}
and integration by parts yields
\begin{multline*}
 \int_{\Omega_h} \partial_x u \dxy = \int_{\Omega_h} \div\begin{pmatrix} u \\ 0 \end{pmatrix} \dxy = \int_{\partial \Omega_h} \begin{pmatrix} u \\ 0 \end{pmatrix} \cdot \nu \ud \mathcal{H}^1\\
 = \int_{{{\{}}y=h(x){{\}}}} u \nu_x\, \mathrm{d}\mathcal{H}^1 = - \int_\R u(x,h(x))h'(x) \dx.
\end{multline*}
Here we used that $\nu = \frac{(-h'(x),1)}{\sqrt{1+h'^2}}$, and so $\nu_x\, \mathrm{d}\mathcal{H}^1 = -h'(x) { \dx}$. 

By the fundamental theorem of calculus, $u(x,h(x)) = x + \int_0^{h(x)} \partial_y u \dy$, and we deduce
\begin{equation*}
 \int_\R u(x,h(x)) h'(x) \dx = \int_\R x h'(x) \dx + \int_{\Omega_h} \partial_y u(x,y) h'(x) \dxy.
\end{equation*}
For the first integral on the right-hand side, integration by parts implies
\begin{equation*}
 \int_\R x h'(x) \dx = - \int_\R h(x) = -V, 
\end{equation*}
and so altogether,
\begin{equation*}
 \int_{\Omega_h} |\nabla u|^2 \dxy = V - \int_{\Omega_h} \partial_y u(x,y) h'(x) \dxy.
\end{equation*}
We add $S_V = \int_\R h'^2 \dx$ to both sides of the equation to obtain
\begin{equation}\label{equalF}
  E_V + S_V = \int_{\Omega_h} |\nabla u|^2 \dxy + S_V = V + S_V - \int_{\Omega_h} \partial_y u(x,y) h'(x) \dxy.
\end{equation}
To show that the energy of $(u,h)$ is not smaller than $V$, it 
{{suffices to prove that}}
\begin{equation*}
 \int_{\Omega_h} \partial_y u(x,y) h'(x) \dxy \le S_V.
\end{equation*}
{{By Lemma~\ref{internvar}, we have either $\partial_yu\equiv 0$, in which case the inequality holds trivially true, or, using H\"older's inequality, }}
\begin{equation}\label{ineqSV}
 \int_{\Omega_h} \partial_y u(x,y) h'(x) \dxy \le \left(\int_{\Omega_h} (\partial_y u)^2(x,y) \dxy\right)^{1/2} \left(\int_{\Omega_h} h'^2(x) \dxy\right)^{1/2} \le \sqrt{\frac{3}{4}}S_V \sqrt{\sup h},
\end{equation}
which 
implies that $E_V(u,h) + S_V(h) \ge V$ for any admissible $(u,h)  \in \mathcal{A}_V$  which satisfies $\sup h \le 4/3$. Note that by \eqref{interpolationineq},
\[
 \sup h \le \left( \frac{9}{16} \right)^{1/3} V^{1/3} S_V^{1/3} \le \left( \frac{9}{16} \right)^{1/3} V^{2/3}.
\]

Hence, the infimum of the energy is exactly $V$ for $0 < V \le \frac{2^5}{9\sqrt{3}}$, and so there exists a maximal $\overline{V} > 0$ such that $F_s(V) = V$ for $V \in (0,\overline{V}]$. 

To show non-existence of a minimizer, let us argue by contradiction and assume that there exists a minimizer $(u,h)  \in \mathcal{A}_V$ for some $V$ {{with}} $0 < V < \overline{V}$. Choose $V_0$ {{with}} $V < V_0 \le \overline{V}${{, and set $\lambda:=\sqrt{\frac{V}{V_0}}$. Then by Lemma \ref{rescale}, since $(u_\lambda,h_\lambda)\in\mathcal{A}_{V_0}$ and $S_V(h) > 0$,}}
we obtain a contradiction:
\begin{eqnarray*}
V_0 &=& F_s(V_0) {{\leq E_{V_0}(u_\lambda,h_\lambda)+ S_{V_0}(h_\lambda)=}} \frac{V_0}{V}E_V(u,h) + \left( \frac{V_0}{V} \right)^{1/2} S_V(h)\\
& <& \frac{V_0}{V} (E_V(u,h) + S_V(h)) = V_0.
\end{eqnarray*}
\end{proof}
We {{next}} consider the regime $\{V>\overline{V}\}$ and aim 
{{to prove}} that minimizers exist {{for every $V>\overline{V}$}}. For that, we need some auxiliary properties. {{For $V>0$, we set
\begin{eqnarray}\label{eq:beta}
\beta(V):=\frac{F_s(V)}{V}. 
 \end{eqnarray}
}}
By \eqref{concavF} and $F_s(V)\le V$, we see that 
$\beta(V)\le 1$ {{for every $V$, and that $\beta$}} is a non-increasing function of V. Using that $F_s(V)=V$ for 
$V \in (0,\overline{V}]$, we can say more about $\beta$:

\begin{Lemma}\label{decrbeta}
 The function $\beta$ is 
strictly
 decreasing in the region $\{F_s(V)<V\}=\{\beta<1\}$.
\end{Lemma}
\begin{proof}
 We assume the contrary, i.e., {{that there exist $\overline{V}<V_0<V_1<\infty$ and $\beta_0<1$ such that}} $F_s(V) = \beta_0 V < V$ for {{all}} $V \in[V_0,V_1] \subset (\overline{V},\infty)$. 
 We use the concavity of $F_s$, $F_s(0) = 0$, and the previous assumption, to get for every $V \in (0,\overline{V})$
 \begin{gather*}
  F_s(V) \ge (1-\frac{V}{V_0})F_s(0) + \frac{V}{V_0} F_s(V_0) = {{\frac{V}{V_0}F_s(V_0)=}}\beta_0 V,\\
\textrm{and}
\\ \beta_0V_0 = F_s(V_0) \ge \frac{V_1-V_0}{V_1 - V} F_s(V) + \frac{V_0-V}{V_1 - V} F_s(V_1).
 \end{gather*}
 {{Since $F_s(V_1)=\beta_0 V_1$, t}}he second relation 
simplifies to $F_s(V) \le \beta_0 V$, which together with the first relation implies $F_s(V) = \beta_0 V$. This is a contradiction, since $\beta_0 < 1$ and by assumption  $F_s(V) = V$.

\end{proof}

The essential step to prove existence of minimizers is to derive compactness of minimizing sequences. It is based on the following auxiliary lemma:
\begin{lemma1}\label{lm:cutting}
 Let $V > \overline{V}$ and $\delta > 0$. Then there exist a length $l=l(V,\delta) > 0$ and {{$C(V,\delta)>0$ 
with the following property: 
For every $(u,h)\in\mathcal{A}_V$ 
with }}$\eps := E_V(u,h) + S_V(u) - F_s(V)\le C(V,\delta)$, {{and every $x_0 < x_1$, with}} $x_1 - x_0 = l$, 
{{we have }}
 \begin{equation}
  \int_{-\infty}^{x_0} h \dx\le \delta \quad \textrm{or\ } \int_{x_1}^\infty h \dx \le \delta.
 \end{equation}
\end{lemma1}

\begin{proof} 
{{Let $V > \overline{V}$ and $\delta > 0$. Note that by the strict monotonicity of $\beta$ in $[\overline{V},\infty)$, we have $\beta(V)<\beta\lt(\frac{V}{2}\rt)$. 
Choose $C(V,\delta)>0$ and let $0<\alpha<1$ be such that
\begin{eqnarray}
 \label{eq:choicepara}
C(V,\delta)<\frac{\delta}{2}\left(\beta\lt(\frac{V}{2}\rt)-\beta(V)\right),\qquad \text{and\quad}\alpha<\left(\frac{\delta}{156}\lt[\beta\lt(\frac{V}{2}\rt)-\beta(V)\rt]\right)^2.
\end{eqnarray}
Assume further that  $\alpha>0$ is such that $\frac{V}{3\alpha}=n\in\N$. Set $l:=\frac{V}{\alpha}$, and let $[x_0,x_1]$ be an arbitrary interval of length $x_1-x_0=l$. Since $\alpha^{-1}V=3n$ with $n\in\N$, we can decompose $[x_0,x_1]$ as an essentially disjoint union of $3n$ intervals of length $1$, i.e.,
\[[x_0,x_1]=\bigcup_{k=0}^{3n-1}[x_0+k,x_0+(k+1)]. \]
}}
Since $\int_\R h \dx = V = n \cdot 3\alpha$, there are at most $n$ intervals which satisfy $\int_{I_k} h \dx \ge 3\alpha$. Similarly, since $S_V(h) < V$, we find at most $n$ intervals such that $\int_{I_k} h'^2 \dx \ge 3\alpha$. Hence, among $I_k$ there is an interval $I$ such that 
\begin{eqnarray}\label{eq:intI}
 \int_I h \dx\le 3\alpha  \qquad \textrm{and\  } \qquad \int_I h'^2\dx\le 3\alpha,
\end{eqnarray}
in particular
 \[3\alpha \ge \int_I h'^2 \dx \ge (\sup_I h -\inf_I h)^2.\]
 Since we also have $3\alpha \ge \int_I h \dx \ge  \inf_Ih$, we get that 
 \begin{equation}\label{suph313}
 \sup_I h = (\sup_I h - \inf_I h) + \inf_I h \le \sqrt{3\alpha} + 3\alpha \le 3(\alpha + \alpha^{1/2}). 
 \end{equation}
{{Without loss of generality (translating $h$), we may assume that $I=[0,1]$. Now we {``cut}'' the profile into three parts. Precisely, we set
\[h_0(x):=\begin{cases}
           h(x),&\text{if\ }x\leq 0,\\
\min(-2h(0)x+h(0),h(x)),&\text{if\ }0\leq x\leq \frac{1}{2},\\
0&\text{otherwise,}
          \end{cases}\text{\qquad and\quad }u_0:=u|_{\Omega_{h_0}},
 \]
and 
\[h_1(x):=\begin{cases}
           \min(2h(1)x-h(1),h(x)),&\text{if\ }\frac{1}{2}\leq x\leq 1,\\
h(x),&\text{if\ }x\geq 1,\\
0&\text{otherwise,}
          \end{cases}\text{\qquad and\quad }u_1:=u|_{\Omega_{h_1}}.
 \]
We use the notation
\[V_0:=\int_{-\infty}^{1/2}h_0\dx,\quad V_1:=\int_{1/2}^\infty h_1\dx,\quad V_{lo}:=V-V_0-V_1.\]
We may assume without loss of generality that $V_0\leq V_1$, which implies that $V_0\leq\frac{V}{2}$.
Note that by construction $(u_0,h_0)\in\mathcal{A}_{V_0}$, $(u_1,h_1)\in\mathcal{A}_{V_1}$, 
\begin{eqnarray}
 \label{eq:Vlo}
V_{lo}\leq \sup_I h\leq 3(\alpha+\alpha^{1/2}),
\end{eqnarray}
and 
\begin{eqnarray}
 \label{eq:E01}
E_V(u,h)\geq E_{V_0}(u_0,h_0)+E_{V_1}(u_1,h_1).
\end{eqnarray}
Further, again by construction and \eqref{suph313},
\begin{eqnarray}\label{eq:S01}
 S_V(u,h)&\geq&S_{V_0}(h_0)+S_{V_1}(h_1)-\left[\int_0^{1/2}h_0'^2\dx+\int_{1/2}^1h_1'^2\dx \right]\nonumber\\
&\geq& S_{V_0}(h_0)+S_{V_1}(h_1)-(2\sup_I h)^2\geq S_{V_0}(h_0)+S_{V_1}(h_1)-36(\alpha+\alpha^{1/2})^2.
\end{eqnarray}
}}

Since $F_s(V)=\beta(V) V$ with $\beta< 1$
 we have by \eqref{eq:E01} and \eqref{eq:S01}, 
\begin{align*}
\beta(V)V + \eps &= F_s(V) + \eps = E_V(u,h) + S_V(h) \ge F_s(V_0) + F_s(V_1) - {{36(\alpha+\alpha^{1/2})^2}}
\\ &= \beta(V_0) V_0 +\beta(V_1) V_1 -{{36(\alpha+\alpha^{1/2})^2}}.
\end{align*}
 {{Moreover, since $\beta$ is non-increasing and since $V=V_0+V_1+V_{lo}$, and $V_0\leq \frac{V}{2}$,}}
\begin{align*}
\beta(V_0) V_0 +\beta(V_1) V_1 &=   \beta(V) V+V_0(\beta(V_0)-\beta(V))+   V_1(\beta(V_1)-\beta(V)) - \beta(V)V_{lo}
\\ &\ge \beta(V) V + V_0(\beta(V/2) - \beta(V)) - V_{lo} 
\end{align*}
which together implies $V_0(\beta(V/2)-\beta(V)) \le \eps + V_{lo} + {{36(\alpha+\alpha^{1/2})^2}}$. Since $V > \overline{V}$, the 
strict
 monotonicity of $\beta$ implies $\beta(V/2) - \beta(V) > 0$, and so, {{since $\alpha^2<\alpha <\alpha^{1/2}$,}}
\begin{equation}\label{v0}
 V_0 \le \frac{\eps + V_{lo} + {{36(\alpha+\alpha^{1/2})^2}}}{\beta(V/2) - \beta(V)} \le 
{{\frac{\eps+78\alpha^{1/2}}{\beta(V/2) - \beta(V)}}}.
\end{equation}
{{By \eqref{eq:choicepara}, we obtain that $V_0<\delta$, and thus
\[\int_{-\infty}^{x_0}h\dx=\int_{-\infty}^{x_0}h_0\dx\leq \int_{\R}h_0\dx=V_0\leq\delta. \]
}}
\end{proof}

\begin{prop1}\label{existminss}
 For any $V > \overline{V}$, a minimizer of \eqref{funcss} exists. 
\end{prop1}
\begin{proof}
Let $V$ be such that $F_s(V)<V$, and let $(u_n,h_n)$ be a minimizing sequence. First we claim that Lemma~\ref{lm:cutting} implies tightness of $(u_n,h_n)$. 

Indeed, let $0 < \delta < V/2$ be fixed and let $l=l(V,\delta)$ be obtained from Lemma~\ref{lm:cutting}. Then for $n$ large enough the energy of $(u_n,h_n)$ is close enough to $F_s(V)$ (i.e., $\eps := E_V(u_n,h_n) + S_V(h_n) - F_{{{s}}}(V)$ from Lemma~\ref{lm:cutting} is small enough). Then we choose $x_{n}$ such that $\int_{-\infty}^{x_{n}} h \dx= 2\delta$ and observe that by Lemma~\ref{lm:cutting}
\begin{equation}\nonumber
 \int_{x_{n}+l}^{\infty} h \dx \le \delta.
\end{equation}
Hence, for $n$ large enough we have $\int_{-\infty}^{x_{n}} h \le 2\delta$ and $\int^{\infty}_{x_{n}+l} h \le \delta$, which implies tightness of the minimizing sequence (up to translations). The existence of a minimizer then follows from the lower semicontinuity of the energy (see Proposition \ref{densityss}).
\end{proof}
We will prove in Corollary \ref{existmincrit} that also for $V=\overline{V}$ a minimizer exists.

\subsection{Regularity of minimizers}\label{sec:reg}

{\bf Notation:} In this section, for given $x \in \R$ ($x_i \in \R$, etc.), we will denote by $z$ ($z_i$) a point in $\R^2$ defined by ${\bf z} := (x,y) = (x,h(x))$ (${\bf z_i} = (x_i,y_i) = (x_i,h(x_i))$).
%

In this subsection we prove the regularity of minimizers of \eqref{funcss} {{if they exist}}. For this we follow the strategy of \cite{FFLM:07} (see also \cite{DeMFus:12}) which in turn was inspired by \cite{ChambLars:03}. Let us first notice that as  in \cite{FFLM:07}, the volume constraint can be relaxed. 
For $\mu>0$ and $V>0$ set 
\[F_s^\mu(V):={{\inf\left\{ \int_{\Omega_h} |\nabla u|^2 \dxy +\int_{\R} h'^2 \dx + \mu\lt|\int_\R h \dx -V\rt|:\ (u,h)\in\mathcal{A}\right\}}}\]
{{with
\[
 \mathcal{A}:=\left\{ (u,h):\ h\in H^1(\R),\ h\geq 0, \ u\in H^1(\Omega_h),
\quad   u(x,0)= x \ \text{if\ }x\in\supp h\right\}. 
\]
Note that $\mathcal{A}_V=\{(u,h)\in\mathcal{A}:\ \int_{\R}h\dx=V\}$.
We have the following relation between $F_s^{\mu}$ and $F_s$:
}}

\begin{prop1}\label{relaxss}
{{Let $C\geq \frac{2}{c_0}$, where $c_0$ is the constant from Proposition \ref{prop:scaling}. Then, }}
for  $\mu := C\min{{\{}}1,V^{-1/5}{{\}}}$, we have $F_s^\mu(V)=F_s(V)$.
\end{prop1}
\begin{proof}
We closely follow the proof of \cite[Theorem 2.8]{golnov}. {{Let $V>0$.
First, since $\mathcal{A}_V\subset\mathcal{A}$, and 
\[\int_{\Omega_h} |\nabla u|^2 \dxy +\int_{\R} h'^2 \dx + \mu\lt|\int_\R h \dx -V\rt|=\int_{\Omega_h} |\nabla u|^2 \dxy +\int_{\R} h'^2 \dx  \]
for every $(u,h)\in\mathcal{A}_V$, we have $F_s^{\mu}(V)\leq F_s(V)$. Assume now for the sake of contradiction that}} $F_s^\mu(V) < F_s(V)$, i.e., there exists $(u,h){{\in\mathcal{A}}}$ such that\footnote{We denote here for simplicity $E(u,h):=\int_{\Omega_h} |\nabla u|^2\dxy$ and $S(h):=\int_{\R} h'^2 \dx$.}
\begin{eqnarray}\label{eq:relaxcontr}
 E(u,h) + S(h) + \mu\lt|\int_\R h \dx - V\rt| < F_s(V).
\end{eqnarray}

{{We note first that b}}y rescaling we can assume that $\int_\R h \dx\le V${{: Indeed, assume that there exists $(u,h)\in\mathcal{A}$ with \eqref{eq:relaxcontr} and $\int_{\R}h\dx=:W>V$. Then for the rescaled pair (see Lemma \ref{rescale}) $(u_\lambda,h_\lambda)\in\mathcal{A}_V$ with $\lambda:=\sqrt{\frac{W}{V}}>1$, we also have by Lemma \ref{rescale} and \eqref{eq:relaxcontr},
\[E(u_\lambda,h_\lambda)+S(h_\lambda)<E(u,h)+S(h)<F_s(V). \]
On the other hand, if $\int_{\R}h\dx\leq \frac{V}{2}$, 
}}
then $\mu |\int_\R h \dx- V|\ge {{\mu\frac{V}{2}=}}\frac{C}{2} \min{{\{}}1, V^{4/5}{{\}}} \ge F_s(V)$ {{by the choice of $C$ and Proposition \ref{prop:scaling}.}}
Hence we can assume that $\int_\R h\dx\ge \frac{V}{2}$. 
{{Now we construct a competitor from $\mathcal{A}_V$ for $(u,h)\in\mathcal{A}\setminus\mathcal{A}_V$}} by rescaling. 
Denoting $\alpha:=\int_{\Omega_h} |\nabla u|^2 \dxy +\int_{\R} h'^2 \dx$ and $W:=\int_\R h \dx$, we find by \eqref{eq:relaxcontr} that
\[\alpha+\mu(V-W) = E(u,h) + S(h) + \mu \lt |\int_{\R} h \dx - V\rt| < F_s(V) \le \alpha \frac{V}{W}.\]
From this {{we get $\mu<\frac{\alpha}{W}$ (since otherwise $\alpha+\mu(V-W)\geq\frac{\alpha V}{W}$)}}, and from the scaling law (Proposition~\ref{prop:scaling}) finally {{(recall that $W\geq V/2$)}}
\[\mu < \frac{\alpha}{W} \le 
 \frac{F_s(V)}{W} {{\leq \frac{1}{W}\frac{1}{c_0} \min\{V,V^{4/5}\}\leq\frac{2}{c_0}\min\{1,\,V^{-1/5}\}  }},\]
 from which we get a contradiction by {{the choice of $C$}}.

\end{proof}

In order to prove more properties of the minimizers, we will need to use the Euler-Lagrange equation. This requires to know some smoothness of the minimizer of~\eqref{funcss}.

\begin{teorema1}\label{thm:reg}
 {{If $(u,h)\in\mathcal{A}$ is a}}  minimizer of \eqref{funcss}{{, then $h$}} is analytic in $\{h>0\}$ and satisfies the zero contact angle condition.
\end{teorema1}

 Following \cite{ChambLars:03,FFLM:07,DeMFus:12}, we divide the proof of the regularity of $(u,h)$ into several lemmas. We first prove that $h$ is locally Lipschitz continuous. Then, in the spirit of \cite{ChambLars:03}, we prove a uniform sphere condition. From this, we derive decay estimates for $|\nabla u|$ which as in \cite{FFLM:07} leads to the regularity of $h$.
 \begin{Lemma}\label{lm:lip}
  Let the pair $(u,h)$ be a minimizer of~\eqref{funcss}. Then in the set $\{ h > 0 \}$ the height profile $h$ is a locally Lipschitz function. 
 \end{Lemma}

 \begin{proof}
 Let $m > 0$, and 
$x_0, x_1 \in \supp h$ be such that $x_0 < x_1$, $ h(x_0) < h(x_1)$ and $m\le h(x)$ for $x\in[x_0,x_1]$. Since $h$ is continuous, there exist {{
\begin{eqnarray*}
 \overline{x}_0&:=&\max\{x\in[x_0,x_1):\ h(x)=h(x_0) \},\\
x_3&:=&\min\{x>x_1:\ h(x)=h(x_0) \},\text{\quad and}\\
x_2&:=&\max\{x<x_3:\ h(x)=h(x_1)\}.
\end{eqnarray*}
Note that }}$h \ge h(x_0)$ in $[\bar x_0,x_3]$ and $h> h(x_0)$ in $(\bar x_0, x_1)$. 
We denote $\delta := h(x_1) - h(x_0).$
 
 Let us now define a competitor $(\tilde u,\tilde h)$ for $(u,h)$. 
In $\R\setminus [\bar x_0, x_3]$ we set $\tilde h:=h$, and in $[\bar x_0, x_3]$
we set
 \begin{equation*}
  \tilde h(x) := \begin{cases} 
h(x_0) & \textrm{ if } h(x_0) \le h(x) < h(x_1) \\ 
h(x) - \delta & \textrm{ if } h(x) \ge h(x_1).
          \end{cases}
 \end{equation*}
 Since $\tilde h \le h$, we have $\Omega_{\tilde h} \subset \Omega_h$ and the following definition makes sense 
 \begin{equation}\label{eq39}
  \tilde u (\x) := u(\x), \quad \x \in \Omega_{\tilde h}.
 \end{equation}

 We denote $M := \{ x \in [\bar x_0,x_3] : h(x_0) \le h(x) < h(x_1) \}$. Then $h'(x) = \tilde h'(x)$ for a.e. $x \not\in M$ and $\tilde h'(x) = 0$ for a.e. $x \in M$, and so
 \begin{equation*}
  \int_\R \tilde h'^2 \dx= \int_\R h'^2 \dx- \int_M h'^2 \dx\le \int_\R h'^2 \dx- \int_{\bar x_0}^{x_1} h'^2 \dx\le \int_\R h'^2 \dx- \frac{\delta^2}{x_1 - \bar x_0},
 \end{equation*}
 where the last inequality follows from  $\delta^2 = (h(x_1) - h(\bar x_0))^2 = \left( \int_{\bar x_0}^{x_1} h' \dx \right)^2 \le (x_1-\bar x_0) \int_{\bar x_0}^{x_1} h'^2 \dx$. \\
 
 We now estimate $\int_\R \tilde h \dx$:
 \begin{equation*}
  \int_\R \tilde h \dx= \int_\R h \dx- \int_{\{ \bar x_0 < x < x_3 \}} \min \{h(x) - h(x_0),\delta\} \dx \ge V - (x_3-\bar x_0)\delta.
 \end{equation*} 
 Since $h(x) \ge m$ in $[\bar x_0, x_3]$, we have that $V \ge \int_{\bar x_0}^{x_3} h(x) \dx \ge (x_3-\bar x_0)m$, and so
 \begin{equation}\label{eq313}
  V \ge \int_\R \tilde h \dx\ge V - \frac{\delta V}{m}.
 \end{equation}
 Finally, by~\eqref{eq39} we see that $\int_{\Omega_h} |\nabla u|^2 \dxy \le \int_{\Omega_{\tilde h}} |\nabla \tilde u|^2 \dxy$.\\
Since $(u,h)$ is a minimizer, the previous estimates and Proposition~\ref{relaxss} imply 
 \begin{align*}
  \int_{\Omega_h} |\nabla u|^2 \dxy+ \int_{\R} h'^2 \dx &\le \int_{\Omega_{\tilde h}} |\nabla \tilde u|^2 \dxy+\int_{\R} \tilde h'^2 \dx+ \mu\lt|\int_\R \tilde h \dx-V\rt| 
\\ &\le \int_{\Omega_h} |\nabla u|^2 \dxy+ \int_{\R} h'^2 \dx- \left( \frac{\delta^2}{x_1-\bar x_0} - \mu \frac{\delta V}{m}\right),
 \end{align*}
 where $\mu = CV^{-1/5}$. Hence $\frac{\delta}{x_1-\bar x_0} \le \mu \frac{V}{m}$, 
 which implies 
 \begin{equation*}
  \frac{h(x_1) - h(x_0)}{x_1 - x_0}\le \frac{h(x_1) - h(x_0)}{x_1 - \bar x_0} = \frac{\delta}{x_1 - \bar x_0} \le \frac{\mu V}{m} \lesssim \frac{V^{4/5}}{m}.
 \end{equation*}

 In the case $h(x_0) > h(x_1) \ge m$ we proceed analogously. Altogether we get that if both $x_0$ and $x_1$ belong to the set $\{ x : h(x) \ge m \}$, then 
 \begin{equation*}
  |h(x_0) - h(x_1)| \lesssim |x_0 - x_1| \frac{V^{4/5}}{m},
 \end{equation*}
 in particular $h$ is locally Lipschitz in the set $\{ h > 0 \}$. 
\end{proof}
We now prove that the graph of $h$ satisfies a uniform sphere condition. The proof is inspired by, but slightly different from, the proof of a similar statement in \cite{ChambLars:03} (see also \cite{FFLM:07, DeMFus:12}). The main difference {{to}} the aforementioned papers is that in our setting, the surface energy is not invariant by rotation of the axis.

\begin{Lemma}\label{lm:ball1}
 Let {{$V>0$, and let}} the pair $(u,h)$ be a minimizer of~\eqref{funcss}. Then, there exists a radius $r_0 = r_0(V) > 0$ {{with the following property: For}} 
 every circle $S_r(\x_0)$ ($\x_0 := (x_0,y_0)$) and every interval $(a,b) \subset \R$ such that $(a,h(a)) \in S_r(\x_0)$, $h(a) > y_0$, $(b,h(b)) \in S_r(\x_0)$, $h(b) > y_0$, and such that the graph of $h$ is above $S_r(\x_0)$ in $(a,b)$, we have that $r > r_0$. 
\end{Lemma}

\begin{proof}
{{Let $S_r(\x_0)$ and $(a,b)\subset\R$ be as in the lemma. }}
We define $(\tilde h,\tilde u)$, a competitor for $(u,h)$, by
\begin{equation*}
 \tilde h(x) := \begin{cases} h(x)& \textrm{ if } x \not\in (a,b)\\ h(a) + \frac{h(b)-h(a)}{b-a}(x-a)& \textrm{ if } x \in (a,b), \end{cases}
\end{equation*}
and $\tilde u(\x) := u(\x)$ for $\x \in \Omega_{\tilde h} \subset \Omega_h$. Since $(u,h)$ is a minimizer, Proposition~\ref{relaxss} implies
\begin{equation*}
 \int_{\Omega_h} |\nabla u|^2 \dxy+ \int_{\R} h'^2 \dx\le \int_{\Omega_{\tilde h}} |\nabla \tilde u|^2 \dxy+ \int_{\R} \tilde h'^2 \dx+ \mu \left| \int_{\R} \tilde h \dx- V\right|.
\end{equation*}
Since $u=\tilde u$ in $\Omega_{\tilde h}$ and $\Omega_{\tilde h} \subset \Omega_h$, we see that $\int_{\Omega_{\tilde h}} |\nabla u|^2 \dxy\le \int_{\Omega} |\nabla \tilde u|^2 \dxy$. We use this together with the fact that $h = \tilde h$ outside of $(a,b)$ to derive
\begin{equation}\label{319}
 \int_a^b \lt(h'^2 - \tilde h'^2\rt) \dx\le \mu \left| \int_{\R} \tilde h \dx- \int_\R h\dx \right| = \mu \left| \int_a^b \lt(h - \tilde h\rt) \dx\right|.
\end{equation}
Using the definition of $\tilde h$ and H\"older's inequality, we obtain for every $x \in (a,b)$
\begin{multline}\label{320}
 \left| h(x) - \tilde h(x) \right| \le \int_a^x \left|h'(x') - \frac{h(b)-h(a)}{b-a}\right| \ud x' \\
\le (b-a)^{1/2} \left( \int_a^b \left( h'(x') - \frac{h(b)-h(a)}{b-a}\right)^2 \ud x' \right)^{1/2}.
\end{multline}
We observe that $\int_a^b h' \dx = h(b) - h(a)$ implies that 
  $$\int_a^b \left( h'(x') - \frac{h(b)-h(a)}{b-a}\right)^2 \ud x' = \int_a^b \ \lt(h'^2(x') - \tilde h'^2(x')\rt) \ud x'.$$ 
 We plug this relation into~\eqref{320} to show
\begin{equation*}
 \int_a^b \lt( h'^2 - \tilde h'^2 \rt)\dx \overset{\eqref{319}}{\le} \mu \int_a^b |\tilde h - h| \dx\le \mu (b-a)^{3/2} \left( \int_a^b \lt(h'^2 - \tilde h'^2\rt) \dx \right)^{1/2}.
\end{equation*}
Hence
\begin{equation*}
 \int_a^b \lt( h'^2 - \tilde h'^2 \rt)\dx \le \mu^2 (b-a)^{3},
\end{equation*}
and subsequently 
\begin{equation}\label{323}
 \int_a^b \lt(h - \tilde h \rt)\dx \le \mu (b-a)^3.
\end{equation}
Let us now estimate $\int_a^b \lt(h - \tilde h\rt) \dx$ using the fact that in the interval $(a,b)$ the height function $h$ is above the circle $S_r(\x_0)$, i.e., that $h(x) \ge y_0 + \sqrt{r^2 - (x-x_0)^2} =: f(x)$ for $x \in (a,b)$. The  trapezoidal rule implies that
\begin{equation*}
 \int_a^b f(x) \dx - (b-a)\frac{f(a)+f(b)}{2} = -\frac{(b-a)^3}{12} f''(\xi) 
\end{equation*}
for some $\xi \in (a,b)$. Using that $\int_a^b \tilde h \dx = (b-a)\frac{h(a)+h(b)}{2}= (b-a)\frac{f(a)+f(b)}{2} $, we get that 
\begin{equation*}
 \int_a^b \lt(h - \tilde h\rt) \dx \ge \int_a^b f \dx - (b-a)\frac{f(a)+f(b)}{2} = -\frac{(b-a)^3}{12} f''(\xi). 
\end{equation*}
Finally, we compute $f''(\xi) = -\frac{r^2}{(r^2 - (\xi-x_0)^2)^{3/2}}$ to show
\begin{equation*}
 \int_a^b \lt(h - \tilde h\rt) \dx \ge -\frac{(b-a)^3}{12} f''(\xi) = \frac{(b-a)^3}{12}\frac{r^2}{(r^2 - (\xi-x_0)^2)^{3/2}} \ge \frac{1}{12r} (b-a)^3.
\end{equation*}
We conclude by combining the previous estimate with~\eqref{323} to get $r \ge \mu/12$.

\end{proof}

Arguing as in \cite{FFLM:07,DeMFus:12}, we obtain the following result:
\begin{Lemma}\label{lm:innerball}
 Let the pair $(u,h)$ be a minimizer of~\eqref{funcss}. Then for every point $x\in\supp h$ there exists a ball $B_{r_0}(x_0,y_0) \subset \Omega_h \cup \{ y \le 0 \}$ (with $r_0$ defined in Lemma~\ref{lm:ball1}) such that $\partial B_{r_0}(x_0,y_0) \cap \textrm{ graph } h = (x,h(x))$. 
\end{Lemma}

%
%
From this and \cite[Lemma 3]{ChambLars:03}, we obtain:
\begin{Corollary}\label{cor:chambolle1}
{{Let the pair $(u,h)$ be a minimizer of~\eqref{funcss}.}}  Let $x_0 \in \R$ be such that $h(x_0) > 0$. Then there exists a neighborhood $U{{\subset\R}}$ of $x_0$ such that $h|_U$ is Lipschitz and admits left and right derivatives at every point of $U$, that are respectively left and right continuous. 
\end{Corollary}

The previous corollary implies that to prove that $h$ is a $ C^1$ function  in the set $\{  h > 0 \}$, it suffices to consider points $x_0 \in \R$, $h(x_0) > 0$, 
{{for which }}$h'_+(x_0) \neq h'_-(x_0)${{, where $h_+'$ and $h_{-}'$ denote the right and left derivatives. F}}ollowing \cite{FFLM:07,DeMFus:12} we call such points {\it corner points}.
 Our aim is to prove that if $(u,h)$ is a minimizer {{of \eqref{funcss}}}, then there are no corner points. In order to show this, we first obtain the following estimate on $|\nabla u|$ (which is also an important ingredient in the proof of \cite{FFLM:07,DeMFus:12}):

\begin{Lemma}\label{lm:gradestimate}
 Let  the pair $(u,h)$ be a minimizer of~\eqref{funcss}, and let $x_0$ be a corner point. Then there exist $\alpha > 1$ and $r_1>0$ such that $\int_{B_\rho({\bf{z_0}}) \cap \Omega_h} |\nabla u|^2 \dxy \le \bar C \rho^{\alpha}$ for all $\rho \in (0,r_1)$, where $\bar C:=r_1^{-\alpha}\int_{\Omega_h} |\nabla u|^2 \dxy$. 
\end{Lemma}

\begin{proof} 
 Since $h(x_0) > 0$, $h$ is locally Lipschitz in the neighborhood of $x_0$ and we can find $\delta,L > 0$ such that 
 \begin{equation}\label{eq317}
  |h(x) - h(x')| \le L|x-x'| \qquad \forall x,x' \in (x_0-\delta,x_0+\delta).
 \end{equation}
 Since by assumption, $x_0$ is a corner point, we can choose $r_1$, $0 < r_1 < \min{{\{}}\delta,h(x_0){{\}}}$, small enough such that for every $\rho \in (0,r_1)$ both the graph of $h|_{(x_0,+\infty)}$ and of $h|_{(-\infty,x_0)}$ intersect $S_\rho({\bf z_0})$ exactly once. For $\rho \in (0,r_1)$, let us denote the arc {{of $S_\rho({\bf z_0})$}}, which connects two intersections of the graph of $h$ with $S_\rho({\bf z_0})$, and which belongs to $\Omega_h$ (i.e., the bottom arc), by $A_\rho$. By virtue of~\eqref{eq317} this arc has length at most 
 \begin{equation}\label{arclength}
  \mathcal{H}^1(A_\rho) \le 2\pi \rho(1-\arctan (1/L)). 
 \end{equation}

 For any $a \in \R$ and any $\rho \in (0,r_1)$, since $\Delta (u-a)=0$ in $B_\rho({\bf z_0}) \cap \Omega_h$ and $\partial_\nu (u-a)=0$ on $\partial \Omega_h \cap B_{\rho}({\bf z_0})$, we  get
 \begin{equation}\label{eq318}
  \int_{B_\rho({\bf z_0}) \cap \Omega_h} |\nabla u|^2 \dxy =  \int_{A_\rho} (u-a) \partial_\nu u \ud \mathcal{H}^1.
 \end{equation}
%
Then Poincar\'e's inequality with the optimal constant implies 
 \begin{equation}\label{eq:poinc}
  \int_{A_\rho} (u-\bar u)^2 \ud \mathcal{H}^1 \le \left(\frac{\mathcal{H}^1(A_\rho)}{\pi}\right)^2 \int_{A_\rho} |\partial_\tau u|^2 \ud \mathcal{H}^1 \overset{\eqref{arclength}}{\le} \left(2\rho(1-\arctan(1/L))\right)^2 \int_{A_\rho} |\partial_\tau u|^2 \ud \mathcal{H}^1,
 \end{equation}
 where $\bar u$ denotes the average of $u$ {{o}}n $A_\rho$ and $\partial_\tau u$ denotes the derivative of $u$ in the tangential direction. Hence, 
{{by}} H\"older's and Young's inequality 
{{we get from}}~\eqref{eq318} {{with $a:=\overline{u}$ and~\eqref{eq:poinc}}} 
 \begin{align*}
  \int_{B_\rho({\bf z_0}) \cap \Omega_h} |\nabla u|^2 \dxy & = \int_{A_\rho} (u-\bar u) \partial_\nu u \ud \mathcal{H}^1
\\&\leq {{\left(\int_{A_\rho} (u-\bar u)^2 \ud \mathcal{H}^1 \right)^{1/2}\left( \int_{A_\rho} |\partial_\nu u|^2 \ud \mathcal{H}^1 \right)^{1/2}  }}
\\ &\le 2\rho(1-\arctan(1/L)) \left( \int_{A_\rho} |\partial_\tau u|^2 \ud \mathcal{H}^1 \right)^{1/2} \left( \int_{A_\rho} |\partial_\nu u|^2 \ud \mathcal{H}^1 \right)^{1/2} 
\\ &\le 2\rho(1-\arctan(1/L)) \frac{1}{2} \int_{A_\rho} \left( |\partial_\tau u|^2 + \left|\partial_\nu u\right|^2 \ud \mathcal{H}^1\right)\\
&= \rho (1-\arctan(1/L)) \int_{A_\rho} |\nabla u|^2 \ud \mathcal{H}^1.
 \end{align*}
We let  $F(\rho) := \int_{B_\rho({\bf z_0}) \cap \Omega_h} |\nabla u|^2 \dxy$ and observe that the last estimate can be rewritten as $F(\rho) \le \rho (1-\arctan(1/L))F'(\rho)$.
 By integrating this inequality, we obtain  for any $\rho \in (0,r_1)$, 
 \begin{equation*}
  F(\rho) \le F(r_1) \left( \frac{\rho}{r_1} \right)^\alpha,
 \end{equation*}
 where $1/\alpha = (1-\arctan(1/L)) < 1$. To conclude we observe that $F(r_1) \le \int_{\Omega_h} |\nabla u|^2 \dxy$ implies 
 \begin{equation*}
  \int_{B_\rho({\bf z_0}) \cap \Omh} |\nabla u|^2 \dxy = F(\rho) \le \rho^\alpha F(r_1)r_1^{-\alpha} \le \bar C\rho^\alpha.
 \end{equation*}
\end{proof}

Following \cite[Th. 3.13]{FFLM:07}, we can now prove that in the set $\{ h > 0 \}$ there are no corner points, and so $h \in C^1(\{ h > 0 \})$:

\begin{Lemma}\label{lm:c1}
 Let the pair $(u,h)$ be a minimizer of~\eqref{funcss}. Then in the set $\{ h > 0 \}$ the height profile $h$ is a $C^1$ function.
\end{Lemma}

\begin{proof}
 To prove the lemma it is enough to show that there are no corner points. Let us argue by contradiction and assume that $x_0$ is a corner point. Then, by Corollary~\ref{cor:chambolle1} and Lemma~\ref{lm:gradestimate} there exist $r_1 > 0$ and $\alpha > 1$ such that for $\rho \in (0,r_1)$ 
 \begin{equation}\label{eq335}
  \int_{B_\rho({\bf z_0}) \cap \Omega_h} |\nabla u|^2 \dxy \le \bar C \rho^{\alpha},
 \end{equation}
 and $h|_{(x_0-r_1,x_0+r_1)}$ is a Lipschitz function with right and left derivatives, which are respectively right and left continuous. Moreover, $x_0$ being a corner point implies $h'_+(x_0) \neq h'_-(x_0)$. 

 First, we observe that by Lemma~\ref{lm:innerball}{{,}}  $h'_-(x_0) < h'_+(x_0)$, and so $\eps := (h'_+(x_0) - h'_-(x_0)) / 4 > 0$. Then (possibly by diminishing $r_1$), we can assume that for every $x \in (x_0,x_0+r_1)$ 
 \begin{equation}\label{eq340}
  \left| \frac{h(x) - h(x_0)}{x - x_0} - h'_+(x_0) \right| \le \eps,
 \end{equation}
 and similarly for $x \in (x_0-r_1,x_0)$ 
 \begin{equation}\label{eq341}
  \left| \frac{h(x) - h(x_0)}{x - x_0} - h'_-(x_0) \right| \le \eps.
 \end{equation} 

 By Lemma~\ref{lm:ball1}, for any $\rho \in (0,r_1)$ there exist unique points $x_l${{, $x_r$}}$\in S_\rho({\bf z_0}) \cap \partial \Omega_h$
{{ with} }$x_l < x_0$ and 
$x_r > x_0$. Using $x_l$ and $x_r$ we define $\tilde h$, a competitor for $h$
{{by}}
 \begin{equation*}
  \tilde h := \begin{cases} h(x) & \textrm{ if } x \not\in [x_l,x_r]\\ a(x) & \textrm{ if } x \in [x_l,x_r],\end{cases}
 \end{equation*}
 where $a(x) := h(x_l) + \frac{x-x_l}{x_r-x_l} (h(x_r)-h(x_l))$ is an affine function which connects $(x_l,h(x_l))$ and $(x_r,h(x_r))$. Using~\eqref{eq335} and the fact that $\partial \Omega_h \cap B_\rho({\bf z_0})$ is Lipschitz, we can extend $u$ to $\Omega_{\tilde h}$ (still denoted by $u$) such that 
\begin{equation}\label{eq335+}
  \int_{B_\rho({\bf z_0}) \cap \Omega_{\tilde h}} |\nabla u|^2 \dxy\le \bar C_1 \rho^{\alpha},
 \end{equation}
 where $\bar C_1$ depends only on $\bar C$ and the Lipschitz constant of $h|_{(x_0-r_1,x_0+r_1)}$.  Since $(u,h)$ is a minimizer of~\eqref{funcss}, Proposition~\ref{relaxss} implies
 \begin{equation}\label{eq338}
  \int_{\R} \lt( h'^2 - \tilde h'^2\rt) \dx \le \int_{\Omega_{\tilde h}\setminus \Omega_h} |\nabla u|^2 \dxy+ \mu \int_{\R} \lt(\tilde h - h\rt) \dx.
 \end{equation}
 
 First we will estimate the left-hand side of~\eqref{eq338} from below. The definition of $\tilde h$ implies
 \begin{equation}\label{329}
  \int_{\R}\lt( h'^2 - \tilde h'^2\rt) \dx= \int_{x_l}^{x_0} h'^2 \dx + \int_{x_0}^{x_r} h'^2 \dx- \int_{x_l}^{x_r} \tilde h'^2 \dx.
 \end{equation} 
 We 
{{set}} $d_l := \frac{h(x_0)-h(x_l)}{x_0 - x_l}$ and $d_r := \frac{h(x_r)-h(x_0)}{x_r - x_0}$. Then
 \begin{equation}\label{eq342}
  \int_{x_l}^{x_0} h'^2 \dx + \int_{x_0}^{x_r} h'^2 \dx \ge d_l^2 (x_0 - x_l) + d_r^2 (x_r - x_0),
 \end{equation}
 and
 \begin{equation}\label{eq343}
  \int_{x_l}^{x_r} \tilde h'^2 \dx = \frac{\left( d_l (x_0 - x_l) + d_r (x_r - x_0) \right)^2}{x_r - x_l}.
 \end{equation}
 
 We plug~\eqref{eq342} and \eqref{eq343} into~\eqref{329} to get
 \begin{equation}\nonumber
  \int_{x_l}^{x_r} \lt(h'^2 - \tilde h'^2\rt) \dx \ge d_l^2(x_0 - x_l) + d_r^2(x_r - x_0) - \frac{\left( d_l (x_l - x_0) + d_r (x_r - x_0) \right)^2}{x_r - x_l}.
 \end{equation}
 A simple algebraic manipulation shows that 
 \begin{equation}\nonumber
  d_l^2(x_0 - x_l) + d_r^2(x_r - x_0) - \frac{\left( d_l (x_l - x_0) + d_r (x_r - x_0) \right)^2}{x_r - x_l} = \frac{(x_0 - x_l)(x_r - x_0)}{x_r - x_l}\left( d_r - d_l \right)^2,
 \end{equation}
 and so the estimates on $d_l$ and $d_r$ (see \eqref{eq340} and \eqref{eq341}) imply that 
 \begin{equation}
  \int_{x_l}^{x_r} \lt(h'^2 - \tilde h'^2\rt) \dx \ge \frac{(x_0 - x_l)(x_r - x_0)}{x_r - x_l}( h'_+(x_0) - h'_-(x_0) - 2\eps)^2 \ge 2\min \{x_0 - x_l,x_r - x_0 \} \eps^2,
 \end{equation}
 where 
{{in the last step}} we used that $2ab/(a+b) \ge \min{{\{}}a,b{{\} }}$ and the definition of $\eps$. Let us now observe that since $h$ is Lipschitz, we have that $x_0 - x_l \ge C \rho$ and $x_r - x_0 \ge C \rho$, where $C$ depends only on the Lipschitz constant of $h{{|_{(x_0-r_1,x_0+r_1)}}}$.
Therefore, for any $\rho \in (0,r_1)$ and the corresponding $x_l,x_r$ we 
{{obtain}} that
 \begin{equation}\label{eq346}
  \int_{x_l}^{x_r} \lt(h'^2 - \tilde h'^2\rt) \dx \ge C \rho, 
 \end{equation}
 where $C > 0$ depends on $h|_{(x_0-r_1,x_0+r_1)}$, but not on $\rho$.  

 Finally, we observe that $|\int_\R \lt(\tilde h - h\rt) \dx| = |\int_{x_l}^{x_r} \lt( \tilde h - h \rt) \dx| \le |B_\rho|$, and so~\eqref{eq338}, \eqref{eq335+}, and~\eqref{eq346} imply that there exists $\alpha > 1$ and constants $C>0,\bar C_1$ such that for every $\rho \in (0,r_1)$ we have $C \rho \le \bar C_1 \rho^{\alpha} + \mu \pi \rho^2$,  which 
{{yields a contradiction}} for sufficiently small $\rho$. 
{{This}} concludes the proof.
\end{proof}

We proved that in the set $\{ h > 0 \}$ the height profile $h$ is a $C^1$ function. 
It then follows that in fact $h$ is more regular. Indeed, for any $x_0$ such that $h(x_0) > 0$ we observe that given $\eps > 0$ there exists $r_1 > 0$ such that for any $\rho \in (0,r_1)$ 
\begin{equation}\label{eq348}
 \mathcal{H}^1(S_\rho({\bf z_0}) \cap \Omega_h) \le (1+\eps) \pi \rho.
\end{equation}
Then we can repeat the proof of Lemma~\ref{lm:gradestimate} while replacing~\eqref{arclength} by \eqref{eq348} to show {{the following result.}}

\begin{Lemma}\label{lm:bettergradestimate}
 Let the pair $(u,h)$ be a minimizer of~\eqref{funcss}, and let $x_0 \in \R$ be such that $h(x_0) > 0$. Then for any $0 < \alpha < 2$ there exists $r_1 > 0$ such that for any $\rho \in (0,r_1)$ 
 \begin{equation}\nonumber
  \int_{B_\rho({\bf z_0}) \cap \Omega_h} |\nabla u|^2 \dxy \le \bar C\rho^{\alpha},
 \end{equation}
 where $\bar C = r_1^{-\alpha}\int_{\Omega_h} |\nabla u|^2 \dxy$. Moreover, given $\alpha < 2$, the corresponding $r_1$ depends (through relation~\eqref{eq348}) only on the modulus of continuity of $h'$ in the neighborhood of $x_0$. 
\end{Lemma}

Similar ideas as in the proof of~Lemma~\ref{lm:c1} then give the following result.
\begin{Proposition}
Let the pair $(u,h)$ be a minimizer of~\eqref{funcss}. Then for every $\beta \in (0,1/2)$ 
the height profile $h|_{\{ h > 0 \}}$ is a $C_{\textrm{loc}}^{1,\beta}$ function. 
\end{Proposition}

{
\providecommand{\U}{\mathcal{U}}

\begin{proof}
Let $\bar x \in \R$, $h(\bar x)>0$, and $\alpha \in (0,2)$ be fixed. Then by Lemma~\ref{lm:bettergradestimate} and the fact that $h \in C^1(\{ h > 0 \})$ there exist $\U$, a neighborhood of $\bar x$, {{a}} radius $r_1 > 0$, and a constant $C_1$ such that for any point $x \in \U$ and any $\rho \in (0,r_1)$ we have
\begin{equation}\label{351}
 \int_{B_\rho({\bf z}) \cap \Omega_h} |\nabla u|^2 \dxy \le C_1 \rho^\alpha.
\end{equation}
Moreover, we can assume that $h'$ is bounded in $\U$.

Let $x_0 \in \U$ be fixed. Using standard extension argument and estimate~\eqref{351} we can extend $u$ 
to $B_\rho({\bf z_0}) \setminus \Omega_h$ (the extension still denoted by $u$) such that $\int_{B_\rho({\bf z_0})} |\nabla u|^2 \dxy\le \bar C_1 \rho^\alpha$, where $\bar C_1$ is independent of $\rho$ and {{the}} choice of $x_0$. 

Now let $x_1 > x_0$ {{be}} such that $|{\bf z_0} - {\bf z_1}| < r_1$. We set $\rho := |{\bf z_0} - {\bf z_1}|$ and define $\tilde h$ by
\begin{equation*}
 \tilde h(x) := \begin{cases} h(x) & \text{{if\ }}x \not \in [x_0,x_1]\\
 a(x) & \text{{if\ }}x \in [x_0,x_1],
                \end{cases}
\end{equation*}
where $a(x) := h(x_0) + \frac{x-x_0}{x_1-x_0} (h(x_1)-h(x_0)$ is an affine function connecting ${\bf z_0}$ and ${\bf z_1}$. We observe that $\Omega_{\tilde h} \subset \Omega_h \cup B_\rho({\bf z_0})$, and so $(\tilde h,u|_{\Omh \cup B_\rho})$ is a well defined competitor for $(u,h)$. Since $(u,h)$ is a minimizer of~\eqref{funcss}, Proposition~\ref{relaxss} implies
\begin{equation*}
 \int_{x_0}^{x_1} \lt(h'^2 - \tilde h'^2\rt) \dx \le \int_{B_\rho({\bf z_0})\setminus \Omega_h} |\nabla u|^2 \dxy + \mu |B_\rho({\bf z_0})| \le \bar C_1 \rho^{\alpha} + \mu \pi \rho^2 \le \bar C \rho^{\alpha}.
\end{equation*}
Since $h'$ is bounded in $\U$, we have that $x_1 - x_0 \ge C^{-1}\rho$, where $C$ does not depend on $\rho$ or $x_1$. Then
\begin{equation*}
 \dashint_{x_0}^{x_1} \left( h'(x) - \dashint_{x_0}^{x_1} h' \right)^2 \dx = \dashint_{x_0}^{x_1} \lt(h'^2 - \tilde h'^2\rt) \dx \le C\rho^{-1} \int_{x_0}^{x_1} \lt(h'^2 - \tilde h'^2\rt) \dx,
\end{equation*}
and so 
\begin{equation}\label{eq355}
 \dashint_{x_0}^{x_1} \left( h'(x) - \dashint_{x_0}^{x_1} h' \right)^2 \dx \le C \rho^{\alpha-1}.
\end{equation}

{{A}} relation 
{{similar to~}}\eqref{eq355} holds also for the choice $x_1 < x_0$, and we can use \cite[Th. 7.51]{AFP:00} to conclude that $h \in C^{1,(\alpha-1)/2}(\overline \U)$. 

\end{proof}
}

We showed that $h \in C_{\textrm{loc}}^{1,\beta}(\{ h > 0 \})$ for any $\beta \in (0,1/2)$, and so $(u,h)$ satisfies all the assumptions of the following theorem {{(see }}\cite[Th.7.49]{AFP:00}{{).}}

\begin{Theorem}
 Let $\Omega$ be an open set in $\R^2$, $g \in L^\infty(\Omega)$ and $u \in H^1(\Omega)$ be a solution of the Neumann problem
\begin{equation*}
\begin{array}{rcll}
  -\Delta u &=&g & \textrm{ in } \Omega\\
 \partial_\nu u &=& 0 & \textrm{ on } S.
\end{array}
\end{equation*}
If $S \subset \partial \Omega$ is a $C^{1,\beta}$ curve relatively open in $\partial \Omega$, $\beta < 1$, then $\nabla u$ has a $C^{0,\beta}$ extension up to $S$.
\end{Theorem}

We apply this theorem to show that for every $s>0$, $|\nabla u|$ is $C^{0,\beta}$ in the neighborhood of $\partial \Omega_h \cap \{ y > s \}$. As a consequence we obtain the following:

\begin{prop1}
 A minimizer $(u,h)$ of \eqref{funcss} satisfies the 
Euler-Lagrange equation
\begin{equation}\label{ELss}
 |\nabla u|^2(x,h(x))-2 h''(x)=\lambda \quad \textrm{ for a.e. } x \in \{ h > 0 \},
\end{equation}
where $\lambda{{=\lambda(V)}}$ is a constant {{that depends on $V$, namely }}the Lagrange multiplier associated 
{{to}} the volume constraint.
\end{prop1}

Having~\eqref{ELss}, a simple bootstrap argument implies that in fact $h \in C^{\infty}(\{ h > 0\})$ and $u \in C^{\infty}(\Omega_h)$ (see \cite{FFLM:07} for more details). If $\kappa$ denotes the mean curvature of $\partial \Omega_h$, observing that \eqref{ELss} can be rewritten as
\[\kappa=\frac{|\nabla u|^2-\lambda}{2} \, \nu_y^3,\] 
we see that  \cite[Th. 3.1]{KLM:05} implies that $h$ is analytic in $\{h>0\}$.
To prove Theorem~\ref{thm:reg} it remains to show that $h$ satisfies the zero contact angle condition. To do so we first derive a lemma analogous to
 Lemma~\ref{lm:gradestimate} which applies to the points of contact with the substrate:

\begin{Lemma}\label{lem:morrey}
 Let the pair $(u,h)$ be a minimizer of~\eqref{funcss}, and let $x_0 \in \supp h$ be such that $h(x_0) = 0$. Then $x_0 \in \partial \supp h$ and there exist $r_0 > 0$ and $\bar C$ such that for every $\rho \in (0,r_0)$
 \begin{equation}\label{343}
  \int_{B_\rho({\bf z_0}) \cap \Omega_h} |\nabla u|^2 \dxy \le \bar C \rho^{4/3}.
 \end{equation}
\end{Lemma}

\begin{proof}
{{First note that by Proposition~\ref{connect}, the set $\{h>0\}$ is connected if $(h,u)$ is a minimizer of~\eqref{funcss}. Therefore, if $h(x_0)=0$ and $x_0\in\supp h$, then $x_0\in\partial \supp h$. Consequently, $h$ }}vanishes 
 in $(-\infty,x_0)$ or in $(x_0,\infty)$. 
 Let us now assume that $h$ vanishes in $(-\infty,x_0)$, the other case being symmetric. Then by Lemma~\ref{lm:innerball} and \cite[Lemma 3]{ChambLars:03} there exists {{a}} radius $r_0 > 0$ such that $h|_{(x_0-r_0,x_0+r_0)}$ has left and right derivatives at every point, that are respectively left and right continuous (but could possibly attain infinite values). Moreover, we can assume (by possibly diminishing $r_0 > 0$) that for every $\rho \in (0,r_0)$, the graph of $h$ intersects $S_\rho((x_0,0))$ in exactly two points $(x_0-\rho,0)$ and $(x_r,h(x_r))$.

 Let us fix $\rho \in (0,r_0)$ and the corresponding $x_r$. Since $u$ minimizes the Dirichlet integral {{in $\Omega_h$ subject to boundary conditions}}
 $u(x,0) = x$, we get that 
 \begin{equation*}
  \int_{B_\rho({\bf z_0}) \cap \Omega_h} \nabla u \cdot \nabla (u-x) \dxy = \int_{S_\rho({\bf z_0}) \cap \Omega_h} (u-x) \partial_\nu u \ud \mathcal{H}^1.
 \end{equation*}
Let $\delta := 1/5$. 
{{By }}H\"older's and Young's inequality the previous relation gives
 \begin{multline}\label{360}
  \int_{B_\rho({\bf z_0}) \cap \Omega_h} |\nabla u|^2 \dxy = \int_{S_\rho({\bf z_0}) \cap \Omega_h} (u-x) \partial_\nu u \ud \mathcal{H}^1 + \int_{B_\rho({\bf z_0}) \cap \Omega_h} \nabla u \cdot 
{{\left(\begin{array}{c}
          1\\0   
            \end{array}\right)
}} \dxy 
\\ \le \left( \int_{S_\rho({\bf z_0}) \cap \Omega_h} \left( \partial_\nu u  \right)^2 \ud \mathcal{H}^1 \right)^{1/2} \left( \int_{S_\rho({\bf z_0}) \cap \Omega_h} (u-x)^2 \ud \mathcal{H}^1 \right)^{1/2} \\
+ \delta \int_{B_\rho({\bf z_0}) \cap \Omega_h} |\nabla u|^2 \dxy + \frac{C}{\delta} |B_\rho({\bf z_0})|.
 \end{multline}
  We estimate the second integral on the right-hand side using Wirtinger's inequality, which states that if $f(0)=0$, then
 \begin{equation*}
  \int_0^l f(t)^2 \ud t \le \left( \frac{2l}{\pi} \right)^2 \int_0^l f'(t)^2 \ud t.
 \end{equation*}
 
  More precisely, we 
{{apply }} Wirtinger's inequality 
{{ to }}$u(x,y)-x$ on $S_\rho({\bf z_0}) \cap \Omega_h$ (observe that $u(x,0) - x = 0$ if $x=0$){{, which,}} together with Young's inequality{{, yields}} 
 \begin{multline*}
  \int_{S_\rho({\bf z_0}) \cap \Omega_h} (u-x)^2 \ud \mathcal{H}^1 \le \rho^2 \int_{S_\rho({\bf z_0}) \cap \Omega_h} |\partial_\tau (u-x)|^2 \ud \mathcal{H}^1\\
 {\le} \rho^2 \int_{S_\rho({\bf z_0}) \cap \Omega_h} (1+\delta) (\partial_\tau u)^2 + (1+\delta^{-1}) \ud \mathcal{H}^1 ,
 \end{multline*}
 where $\tau$ denotes the tangent vector to $S_\rho$, and we used that $|\partial_\tau x| \le 1$. Hence, using Young's inequality again and $|S_\rho({\bf z_0}) \cap \Omega_h| \le \rho \pi/2$,
 it follows from~\eqref{360}: 
 \begin{multline*}
  (1-\delta) \int_{B_\rho({\bf z_0}) \cap \Omega_h} |\nabla u|^2 \dxy \le \frac{C}{\delta} \rho^2  + \frac{\rho}{2} \int_{S_\rho({\bf z_0}) \cap \Omega_h}  \left[\left( \partial_\nu u \right)^2 + (1+\delta) (\partial_\tau u)^2 + (1+\delta^{-1}) \right] \ud \mathcal{H}^1
\\ \le C\rho^2 \delta^{-1} + \frac{\rho(1+\delta)}{2} \int_{S_\rho({\bf z_0}) \cap \Omega_h} |\nabla u|^2 \ud \mathcal{H}^1.
 \end{multline*}
 Since $\delta = 1/5$, one has $(1-\delta)/(1+\delta)=2/3$, and so
  \begin{equation*}
  \int_{B_\rho({\bf z_0}) \cap \Omega_h} |\nabla u|^2 \dxy \le C_1\rho^2 + \frac{3}{4} \rho \int_{S_\rho({\bf z_0}) \cap \Omega_h} |\nabla u|^2 \ud \mathcal{H}^1.
 \end{equation*}
 If we denote $G(\rho) := \int_{B_\rho({\bf z_0}) \cap \Omega_h} |\nabla u|^2 \dxy + 2C_1 \rho^2$, then the last relation is equivalent to
 \begin{equation*}
  G(\rho) \le \frac{3}{4} \rho G'(\rho). 
 \end{equation*}
{{Integrating this }}relation (the same way as we did in the end of the proof of Lemma~\ref{lm:gradestimate}) 
{{we}} obtain~\eqref{343}.
\end{proof}

{{Proceeding along the lines of}} the proof of Lemma~\ref{lm:c1}, 
{{we derive the following result from Lemma \ref{lem:morrey}.}}
\begin{Lemma}\label{lm:zerocontact}
{{Let the pair $(u,h)$ be a minimizer of~\eqref{funcss}.}} 
 Let $x_0 \in \supp h$ be such that $h(x_0) = 0$. Then $h'(x_0) = 0$, and so $h \in C^1(\R)$. 
\end{Lemma}

\subsection{More qualitative results}\label{sec:morequal}
It has been observed for various variational models for the epitaxial growth that minimizers are often not unique. 
For instance, for a model for periodic island formation, there is a regime of volumes and periods in which the flat configuration is not the only minimizing 
configuration (see \cite[Theorem 2.13]{FM:12}). We do not obtain uniqueness {{of minimizers}} here, but part of the following result is the weaker statement that for almost every $V$,
 any two minimizers have the same surface and the same elastic energy. For a similar result for faceted islands see \cite{FonPraZwi:14}.

\begin{prop1}
\label{prop:diff}
 {{The function }}$F_s${{$:V\mapsto F_s(V)$}} is  Lipschitz continuous with Lipschitz constant less than $C\min{{\{}}1,V^{-1/5}{{\}}}$. Moreover, if $F_s$ is differentiable at $V$, then 
\begin{eqnarray}\label{eq:Fs'}
 F_s'(V)=\lambda_V=\frac{1}{V}(E_V+1/2 S_V).
\end{eqnarray}
 At such points of differentiability, if $(u,h)$ and $(\tilde u, \tilde h)$ are two different minimizers of \eqref{funcss}, $E_V(u,h)=E_V(\tilde u, \tilde h)$ 
and $S_V(u,h)=S_V(\tilde u, \tilde h)$. Finally, there holds,
\[\varlimsup_{\eps \to 0^+} \frac{F_s(V+\eps)+F_s(V-\eps)-2F_s(V)}{\eps^2} \le -\frac{S_V}{4V^2}.\]
\end{prop1}

\begin{proof}
We already know from Proposition \ref{concavss} that $F_s$ is locally Lipschitz continuous. The estimate on the Lipschitz constant can be obtained by two different 
ways. The first approach is to use Proposition \ref{relaxss}, and test for two volumes $V$ and $W$ the minimization problem $F_s^\mu$ with the minimizers for each of these volumes. Another approach is to compute the derivative of $F_s$ directly. For this, we see that for $\eps>0$, using the rescaling argument (see Lemma \ref{lem:resc}),
\[F_s(V+\eps)-F_s(V)\le (1+\eps/V) E_V+(1+\eps/V)^{1/2}S_V -E_V-S_V,\]
where $E_V$ and $S_V$ are the elastic and surface energy of the minimizer, respectively. If $F_s$ is differentiable at $V$, this implies {{by Proposition~\ref{prop:scaling} that}}
 $F_s'(V)\le \frac{1}{V}(E_V+1/2 S_V)\le F_s(V)/V\le {{c_0^{-1} }} \min{{\{}}1,V^{-1/5}{{\}}}$ {{with $c_0>0$ from Proposition~\ref{prop:scaling}}}.
Similarly, by rescaling from $V$ to $V-\eps$, we find {{by Lemma~\ref{rescale}}} that
\[\frac{F_s(V-\eps)-F_s(V)}{-\eps}\ge \frac{1}{\eps}(E_V+S_V-(1-\eps/V)E_V-(1-\eps/V)^{1/2}S_V),\]
and thus $F_s'(V)\ge \frac{E_V+1/2 S_V}{V}\ge {{\frac{c_0}{2}}}\min{{\{}}1,V^{-1/5}{{\}}}$, which implies $F_s'(V)=\frac{1}{V}(E_V+1/2 S_V)$. Moreover, this also implies that two minimizers for volume $V$ have the same elastic and the same surface energy. The same rescaling argument also gives the bound
\[\varlimsup_{\eps \to 0^+} \frac{F_s(V+\eps)+F_s(V-\eps)-2F_s(V)}{\eps^2} \le -\frac{1}{4V^2}.\]

It remains to show that 
\begin{equation}\label{equallambda}\lambda_V V=F_s(V)-\frac{1}{2} S_V\end{equation}
(where $\lambda_V$ is the Lagrange multiplier)
for 
 $V > \overline{V}$. We test the Laplace equation for $u$ {{in $\Omega_h$}} with the function $y\partial_y u$ to find
\[0=\int_{\Omh} -\div \begin{pmatrix} \partial_x u \\ \partial_y u \end{pmatrix} (y \partial_y u) \dxy=\int_{\Omh} \begin{pmatrix} \partial_x u \\ \partial_y u \end{pmatrix} \cdot \begin{pmatrix}y  \partial_{xy} u \\ y\partial_{yy}u+\partial_y u \end{pmatrix} \dxy. \]
Using integration by parts we obtain
\begin{align*}
 2\int_{\Omh} \partial_x u (y \partial_{xy} u) \dxy&=-\int_{\Omh} (\partial_x u)^2 \dxy+\int_{\partial \Omh} (\partial_x u)^2 y \nu_y \dxy, \text{{\quad and}}
\\ 2\int_{\Omh}  \partial_y u (y\partial_{yy} u) \dxy&=-\int_{\Omh} (\partial_y u)^2 \dxy+\int_{\partial \Omh} (\partial_y u)^2 y \nu_y \dxy,
\end{align*}
which together with the previous relation imply
\[ \int_{\Omh} \lt[(\partial_x u)^2 - (\partial_y u)^2 \rt]\,  \dxy= \int_{\partial \Omh} |\nabla u|^2 y\nu_y \, \ud \mathcal{H}^1. \]
Using the Euler-Lagrange equation, the right-hand side can be written as $-2 \int_\R h'^2 \dx + \lambda \int_\R h \dx$, which implies
\[\int_{\Omega}\lt[(\partial_x u)^2-(\partial_y u)^2\rt] \,  \dxy=-2S_V+\lambda_V V.\]
From this and Lemma \ref{internvar}, we finally obtain \eqref{equallambda}. From \eqref{equallambda} we get that if $F$ is differentiable, then $F'(V)=\lambda_V$.


\end{proof}



We can now use this information to study the compactness properties of minimizers.
\begin{prop1}\label{prop:suppest}
 Let $V > \overline{V}$. Then for every minimizer $(u,h)$ of \eqref{funcss}, the height function $h$ has bounded support, and $\mathcal{H}^1(\supp (h))\le\frac{\lambda_V S_V}{1-\lambda_V}\les V^{3/5}$.
\end{prop1}
\begin{proof}
Let us first prove that any minimizer of~\eqref{funcss} is compactly supported. For the sake of contradiction, assume it is not. Since
\begin{equation}\nonumber
 \int_{\R} h'^2 \dx + \int_\R \lt( \int_0^{h(x)} |\nabla u|^2(x,y) \dy \rt) \dx <+\infty,
\end{equation}
for any $\eps>0$ and $K > 0$ there exist $x_1 <  x_2$ such that $x_2-x_1\ge K$ and
\begin{equation}|h'(x_i)| +\int_0^{h(x_i)} |\nabla u|^2(x_i,y) dy < \eps{{\quad i=1,2}}.
\label{bdsp1}\end{equation}
Using 
\[
u(x_2,h(x_2))-u(x_1,h(x_1))=\int_{x_1}^{x_2} \nabla u (x,h(x))\cdot 
\begin{pmatrix} 1 \\ h'(x) \end{pmatrix}  \dx\]
 we find
\begin{align*}
|u(x_1,h(x_1))-u(x_2,h(x_2))|&\le\int_{x_1}^{x_2} |\nabla u (x,h(x))| \sqrt{1+h'^2(x)} \dx
\\ &\overset{\textrm{H\"older}}{\le} \left(\int_{x_1}^{x_2} |\nabla u|^2(x,h(x)) \dx\right)^{1/2} \left(\int_{x_1}^{x_2} (1+h'^2(x))\dx \right)^{1/2}.\end{align*}
We use the Euler-Lagrange equation{{~\eqref{ELss}}} to replace the first 
{{term}} on the right-hand side 
{{of}} the previous relation, and get 
\begin{align} \label{bdsp2}
|u(x_1,h(x_1))-u(x_2,h(x_2))|&{{\leq \left(\int_{x_1}^{x_2} (\lambda+2h^{\prime\prime}(x)) \dx\right)^{1/2} \left(\int_{x_1}^{x_2} (1+h'^2(x))\dx \right)^{1/2} }}
\nonumber\\
&\le(\lambda (x_2-x_1)+2 |h'(x_1)-h'(x_2)| )^{1/2}( x_2-x_1 +S_V)^{1/2}\nonumber\\ 
&\overset{\eqref{bdsp1}}{\le} (\lambda (x_2-x_1)+2 \eps)^{1/2}( x_2-x_1 +S_V)^{1/2}.\end{align}
Finally, since $u(x,h(x))= x+\int_{0}^{h(x)} \partial_y u(x,y) dy$, we have
\begin{align*}
|u(x_1,h(x_1))-u(x_2,h(x_2))|&\ge |x_1-x_2|-\int_0^{h(x_1)} |\partial_y u(x_1,y)| dy-\int_0^{h(x_2)} |\partial_y u(x_2,y)| dy\\
& \overset{\eqref{bdsp1}}{\ge} |x_1-x_2|-2\eps^{1/2} (\sup h)^{1/2},\end{align*}
so that 
\[x_2-x_1 \le (\lambda (x_2-x_1) +2 \eps)^{1/2}( x_1-x_2 +S_V)^{1/2} + 2\eps^{1/2} \sup h^{1/2}.\]
Since $\lambda < 1$ and $x_2-x_1 \ge K$, we get a contradiction for $K$ large enough. 
We thus see that the support of $h$ must be bounded.
 
{{To}} prove the estimate on the size of the support of $h$, 
 take $x_1 < x_2$  on the boundary of the support. From the zero contact angle condition (see Lemma~\ref{lm:zerocontact}) we get $h'(x_i)=0$. Since $u(x_i,h(x_i))=x_i$, 
~\eqref{bdsp2} {{with $\eps=0$}} implies $x_2-x_1 \le \lambda(x_2-x_1+S_V)$. {{Note that by Proposition~\ref{prop:diff} and the scaling law, $\lambda_V\leq\frac{V^{4/5}}{c_0 V}$. Putting things together, we get}}
\[x_2-x_1 \le \frac{\lambda S_V}{1-\lambda} \lesssim V^{3/5}. \]
\end{proof}
\begin{remark1}
 The 
bound on the size of support of $h$ {{derived in Proposition~\ref{prop:suppest}}} is slightly suboptimal since we expect {{from the proof of the scaling law (see Proposition\ref{prop:scaling})}} that $\mathcal{H}^1(\supp (h))\sim V^{2/5}$.
\end{remark1}

Using 
Proposition{{~\ref{prop:suppest}}}, we can prove existence of a minimizer at the critical volume.
\begin{coroll1}\label{existmincrit}
 There exists a minimizer of \eqref{funcss} for $V=\overline{V}$.
\end{coroll1}
\begin{proof}
 The existence of a minimizer will follow from a general estimate on $\lambda_V$. From \eqref{equallambda} we know that for $V > \overline{V}$ we have $V > F_s(V) = \lambda_VV + S_V/2$, from where we get
\begin{equation}\label{eq:lambdaSV}
 \lambda_V < 1 - \frac{1}{2}\frac{S_V}{V}. 
\end{equation}
To show that $\lambda_V$ 
{{is bounded}} away from $1$ as $V \to \overline{V}$, it is enough to show that $S_V/V$ does not 
{{tend}} to zero as $V \to \overline{V}$. Let us fix $V > \overline{V}$. Then, from \eqref{equalF}, \eqref{ineqSV} {{and $F_s(V)\leq V$}} we get that {{$V\geq V+S_V-\sqrt{\frac{3}{4}}S_V\sqrt{\sup h}$, which implies that}} $\sup h > \frac{4}{3}$. By 
{{Remark \ref{rem:interpol}}} we have $\sup h \le (9/16)^{1/3} V^{1/3} S_V^{1/3}$, {{and}} thus
\begin{equation*}
 \frac{S_V}{V} > \left(\frac{4}{3}\right)^3 \frac{16}{9} V^{-2}.
\end{equation*}
We thus get {{with \eqref{eq:lambdaSV}}}
\begin{equation*}
 \lambda_V \le {{\frac{16}{9}\frac{(\sup h)^2}{V^2}\le }}1 - \frac{2^9}{3^5} {V}^{-2}
\end{equation*}
for any $V > \overline{V}$. 
{{By P}}roposition{{~\ref{prop:suppest},}} 
the size of the support of {{$h$ for }}a{{ny}} minimizer $h$ (for any $V > \overline{V}$) is bounded by $\frac{\lambda_V}{1-\lambda_V} S_V$. 
Hence, the limit of minimizers as $V \to \overline{V}$ exists. This limit is a minimizer of $F_s(\overline{V})$ (as a consequence 
of a  simple $\Gamma-$limit type argument 
along the lines of Proposition \ref{densityss}). 

\end{proof}

 \subsection{Asymptotic analysis}\label{sec:asympss}
For every $V > \overline{V}$ and every 
$(u,h){{\in\mathcal{A}_V}}$ 
with $u$  the minimizer of the Dirichlet energy {{in $\Omega_h$ subject to the boundary condition}}, we define an (anisotropically) rescaled height profile $\widetilde h(x):= V^{-3/5}h(V^{2/5}x)$ and a rescaled energy
\[G_V(\widetilde h):=V^{-4/5}(S_V(h) + E_V(u,h))= \int_{\R} \widetilde h'^2 \dx + V^{-4/5} \int_{\Omh} |\nabla u|^2 \dxy.\]
Observe that the rescaled $\widetilde h$ satisfies $\int_\R \widetilde h\dx =1 $.

\begin{teorema1}\label{th:asympss}
 For every sequence $V_n \to +\infty$ and every minimizer $(h_{V_n},u_{V_n})$ of \eqref{funcss}, the corresponding $\widetilde h_{V_n}$ (possibly translated) converge, up to a subsequence, in $L^\infty(\R)$ to some function $h$, which minimizes the functional
\begin{equation}\label{Gss}
G(h):=\left(\inf_{u(x,0)=x}\int_{\{h>0\}\times[0,+\infty)} |\nabla u|^2 \dxy \right)+\int_\R h'^2 \dx
\end{equation}
 under the constraint $\int_\R h \dx=1$.

\end{teorema1}
\begin{proof}
We start by noticing that Proposition~\ref{prop:scaling} and~\eqref{interpolationineq} give $V^{3/5} \les \max h \les V^{1/3}S_V^{1/3}$, which implies $S_V \gtrsim V^{4/5}$. If $F_s(V)$ is differentiable at $V$, then Proposition~\ref{prop:diff} gives $\beta'(V) = \left( \frac{F_s(V)}{V} \right)' = -\frac{S_V}{2V^2}$. This together with $S_V \gtrsim V^{4/5}$ implies $\beta'(V) \le -CV^{-6/5}$, and so
\begin{equation}\label{estimbeta}
 \beta(V)\le\beta(V/2)-C\int_{V/2}^V t^{-6/5} \ud t=\beta(V/2)-CV^{-1/5}(2^{1/5}-1).
\end{equation}

Since $G_{V_n}(\widetilde h_{V_n}) = V_n^{-4/5}(S_{V_n}(h_{V_n}) + E_{V_n}(u_{V_n},h_{V_n}))$ and by the scaling law (see Proposition~\ref{prop:scaling}) $F_s(V_n) = S_{V_n}(h_{V_n}) + E_{V_n}(u_{V_n},h_{V_n}) \les V_n^{4/5}$, we see that $G_{V_n}(\widetilde h_{V_n}) \le C$. 

Let us prove that (after possible translation) the sequence $\tilde h_{V_n}$ is tight. For this we follow the 
argument 
{{from}} the proof of Lemma \ref{lm:cutting}. Since $\int_\R \widetilde h_{V_n} \dx = 1$ and $\int_\R \widetilde h_{V_n}'^2 \dx \le C$ (independently of $V_n$), for every $\eps>0$ there exists $l = l(\eps)$ such that in every interval of length at least $l$ we can cut the profile $\tilde h_V$ in two parts of volume $\alpha_1$ and $\alpha_2$, respectively, with $1 - (\alpha_1+\alpha_2) \le \eps${{ (for the precise construction we refer to the proof of Lemma~\ref{lm:cutting})}}. Moreover, we can assume that the cost (surface energy) of the cut is bounded by a multiple of $\eps$. Then we get
\[\beta({V_n}){V_n}^{1/5} + C\eps = G_{V_n}(\tilde h_{V_n}) + C\eps \ge {V_n}^{1/5}(\alpha_1\beta({V_n}\alpha_1) + \alpha_2\beta({V_n}\alpha_2)),\]
{{and }}thus
\begin{multline*}\eps(C+V_n^{1/5}\beta({V_n}))\ge V_n^{1/5}(\alpha_1(\beta(V_n\alpha_1)-\beta(V_n))+\alpha_2(\beta(V_n\alpha_2)-\beta(V_n)))\\
 \ge V_n^{1/5}(\alpha_1(\beta(V_n/2)-\beta(V_n)),
\end{multline*}
where we used that $\beta$ is non-increasing and $\alpha_1 < 1/2$. 
Using that $V_n^{1/5}\beta(V_n)\le C$ and \eqref{estimbeta} we find
\[\alpha_1 \le C\eps,\] 
which by the same argument as in the proof Proposition~\ref{existminss} implies tightness of (possibly shifted) $\widetilde h_{V_n}$. 
As a consequence of tightness we get $L^1(\R)$ convergence of (possibly a subsequence of) $\tilde h_{V_n}$ to $h$. Using the compact embedding (on bounded domains) of $H^1$ into $L^\infty$ we get locally uniform convergence of $\widetilde h_{V_n}$ to $h$. Moreover, tightness of the sequence $\widetilde h_{V_n}$ implies that outside of a compact set we can use {{a}} half-line version of \eqref{interpolationineq} to show that $\widetilde h_{V_n}$ is uniformly small there. These two facts together 
{{yield}} uniform convergence.

 Let us now show that 
\[\varliminf_{V_n\to +\infty} G_{V_n}(\widetilde h_{V_n})\ge G(h).\] 
Since the surface energy is lower semicontinuous, to prove the previous relation it is enough to prove the inequality for the elastic part of the energy. For every $\eps>0$ and for $V_n$ large enough we can assume that $\{h>\eps\}\subset \{\tilde h_{V_n}> \eps/2\}$, and so
\[\int_{\Omega_{ h_{V_n}}} |\nabla u_{V_n}|^2 \dxy\ge  \int_{\Omega_{ h_{V_n}}\cap\left[ \{h_{V_n}>\eps V_n^{3/5}/2\}\times[0,\infty)\right]} |\nabla u_{V_n}|^2 \dxy\ge \int_{ \{h>\eps V_n^{3/5}\}\times[0, \eps V_n^{3/5}]} |\nabla u_{V_n}|^2 \dxy.\]
Using the change of variables $x= V_n^{2/5}\hat x$, $ y=V_n^{2/5} \hat y$, and $u_{V_n}=V_n^{-2/5} \hat u_{V_n}$ we find
\[ V_n^{-4/5}\int_{\Omega_{ h_{V_n}}} |\nabla u_{V_n}|^2 \dxy\ge \int_{ \{h>\eps/2\}\times[0, \eps V_n^{1/5}]} |\nabla \hat u_{V_n}|^2 \dxy.\]
Since for any interval $I \subset \R$
\[\lim_{L \to \infty} \min_{u(x,0)=x} \int_{I\times[0,L]} |\nabla u|^2 \dxy = \min_{u(x,0)=x} \int_{I\times[0,\infty)} |\nabla u|^2 \dxy,\]
we obtain 
\[\varliminf_{n\to+\infty} V_n^{-4/5}\int_{\Omega_{ h_{V_n}}} |\nabla u_{V_n}|^2 \dxy \ge \min_{u(x,0)=x} \int_{ \{h>\eps/2\}\times[0, +\infty)} |\nabla  u|^2 \dxy.\]
By letting $\eps\to 0$ we obtain the desired lower bound. \\
For any other admissible function $g$, {{ it is easily seen that 
\[\varlimsup_{V_n\to +\infty} G_{V_n}( g)\le G(g).\]}}
By a classical argument of $\Gamma-$convergence {{(see }}\cite{Br:02}{{)}}, we deduce that $h$ is a minimizer of $G$.
\end{proof}

\begin{remark1}
Notice that
\[G(h):= C_W \sum_{i\in \N} (b_i-a_i)^2+\int_\R h'^2 \dx,\]
where the intervals $(a_i,b_i)$ are the connected components of $\{h>0\}$, 
and $C_W$ is a constant defined by
\begin{eqnarray}\label{eq:CW}
C_W:=\inf \left\{ \int_{[0,1]\times[0,+\infty)} |\nabla u|^2 \dxy : u(x,0)= x \right\}.
\end{eqnarray}
\end{remark1}
{{We}} now study the limiting problem{{.}}
\begin{prop1}
 The minimization problem \eqref{Gss} admits (up to translations) a unique minimizer $\overline{h}$
{{given by}}
\[\overline{h}(x):=\begin{cases}
                    \frac{3}{2}\ell^{-3}  \left(\ell^2-x^2\right){{,}} \quad &\text{{if\ }}x\in[-\ell,\ell]\\[8pt]
		    0{{,}} \qquad \qquad &\text{{if\ }}x\notin[-\ell,\ell],
                   \end{cases}
\]
where $\ell:= \left(\frac{9}{16 C_W}\right)^{1/5}$.
\end{prop1}

\begin{proof}
The existence of a minimizer follows either by 
{{T}}heorem{{~\ref{th:asympss}}} or by the following direct argument, which is similar in spirit to \cite[Proposition 4.5]{GolZwick}. Let $h$ be an admissible function. If $S_i$ denotes the surface energy of the $i-$th connected component of $h$, $\ell_i$ its length, and $V_i$ its volume (assuming $V_i$ are in a non-increasing order), we observe that $h(x) = \int_{a_i}^x h' \dx$ for $x \in [a_i,b_i]$, and so $V_i \le \ell_i^{3/2} S_i^{1/2}$. We sum this inequality for $i \ge n$ and apply H\"older's inequality to get
\[\sum_{i \ge n} V_i\le \Big(\sum_{i \ge n} \ell_i^3\Big)^{1/2} \Big(\sum_{i \ge n} S_i\Big)^{1/2}.\]
Since $\sum_{i \in \N} C_W \ell_i^2 + S_i = G(h) \le C$, we get that $n \ell_n^2 \le C$, i.e., $\ell_n \le C n^{-1/2}$. 
{{We deduce }}$\sum_{i\ge n} V_i \le Cn^{-1/2}$, which shows tightness of a minimizing sequence. 

Let now $\overline{h}$ be a minimizer. Then in each of its connected component $[c_i-\ell_i,c_i+\ell_i]$, 
$\overline{h}$ satisfies $\overline{h}''=-\lambda_i$, and so $\overline{h}=-\frac{\lambda_i}{2}((x-c_i)^2-\ell_i^2)$. 
Since $\int_{c_i-\ell_i}^{c_i+\ell_i} \overline{h} \dx=V_i$, the form of $\overline{h}$ implies $\lambda_i=\frac{3 V_i}{2 \ell_i^3}$. Then by direct computation the energy inside $[c_i-\ell_i,c_i+\ell_i]$ 
{{equals}}
\[4C_W \ell_i^2+\frac{2}{3}\lambda_i^2 \ell_i^3=4C_W \ell_i^2+\frac{3 V_i^2}{2 \ell_i^3}.\]
This expression is minimized (under the constraint of volume $V_i$) by $\ell_i=\left(\frac{9}{16C_W}\right)^{1/5} V_i^{2/5}$. Then the total energy $G(\overline{h})$ is 
\[\sum_i \left( 324 \cdot C_W^3 \right)^{1/5} V_i^{4/5}\]
with the constraint $\sum V_i=1$. Thus, 
this energy is minimized by a single island.\\
\end{proof}
\begin{remark1}\rm
 By the uniqueness of the minimizer of $G$, we see that the whole sequence $\widetilde h_V$ (possibly translated) converges in $L^\infty$ to $\bar h$.
\end{remark1}
{{We will finally}} prove 
the 
exponential {{rate of }}convergence {{of optimal profiles. For that, we need}}, 
the following quantitative inequality, which can be considered as a very simple quantitative isoperimetric inequality {{(see }}\cite{FMP:08,CFMP:09}{{)}}.
\begin{prop1}
\label{prop:exp1}
 Let $L>0$, $V>0$, and let $h_{\min}{{\in H^1(\R)}}$ be the minimizer of $\int_{-L}^L h'^2 \dx$ under the constraints $h(-L)=h(L)=0$ and $\int_{-L}^L h \dx =V$. Then for every other $h{{\in H^1(\R)}}$ satisfying the same constraints,
\begin{equation}\label{isoper}
 \int_{-L}^L h'^2 \dx-\int_{-L}^L h_{\min}'^2 \dx\ge \frac{1}{4L^2} \int_{-L}^L|h-h_{\min}|^2 \dx.
\end{equation}
\end{prop1}

\begin{proof}
Let $L,V,$ and $h_{min}$ be as in the statement. Then, for every competitor $h$ we write 
$h=(h-h_{\min})+h_{\min}$, and so 
$\int_{-L}^L h'^2 \dx= \int_{-L}^L (h-h_{\min})'^2 \dx+\int_{-L}^L h_{\min}'^2 \dx+2\int_{-L}^L h_{\min}'(h-h_{\min})' \dx$. Since $h_{min}''$ is constant and $\int_{-L}^L h \dx= \int_{-L}^L h_{min} \dx$, integration by parts implies
 \[\int_{-L}^L h'^2 \dx= \int_{-L}^L (h-h_{\min})'^2 \dx +\int_{-L}^L h_{\min}'^2 \dx.\]
For $x\in [-L,L]$ we have
\[|h(x)-h_{\min}(x)|^2\le \left(\int_{-L}^x|h'-h_{\min}'| \dx\right)^2\le 2L \int_{-L}^L (h-h_{\min})'^2 \dx,\]
and so by integration we obtain
\[ \int_{-L}^L |h - h_{min}|^2 \dx \le 4L^2 \int_{-L}^L (h-h_{min})'^2 \dx = 4L^2 \left( \int_{-L}^L h'^2 \dx-\int_{-L}^L h_{\min}'^2 \dx \right) ,\]
{{which shows the claim.}}
\end{proof}

We now prove the exponential convergence of $\widetilde h_V$ to a truncated parabola. To state our result, we will need the following notation. Let $V>\overline{V}$ and $\widetilde h_{V}$, a  minimizer of $G_{V}$, be fixed. Then for $s>0$, we let $\widetilde I_{s}$ be the largest connected component of $\{\widetilde h_{V}>s\}$ and $\bar h_{s}$ be the minimizer of $\int_{\widetilde I_{s}} h'^2 \dx$ with the constraint $\int_{\widetilde I_{s}} h \dx=\int_{\widetilde I_{s}} \widetilde h_V \dx$ and $h=\widetilde h_V$ on the boundary of $\widetilde I_s$.
{\newcommand{\hh}{\widetilde h_V}

\begin{prop1}\label{expss}
 Let $\eps > 0$. Then there exist constants $C_0=C_0(\eps)$ and $C_1=C_1(\eps)$ such that for every $V > \overline{V}$ and for every minimizer $\hh$ of $G_V$,
\[\|\hh-\bar h_s\|_{L^2(\widetilde I_s)}\le C_0 \exp(- C_1 V^{1/5})\quad \forall s \ge \eps.\] 
\end{prop1}
\begin{proof}
Let $\eps>0$, $V > \overline{V}$, and $\hh$ be as in the statement. If $V$ is large enough, then 
\begin{equation*}
 \|\hh - \overline{h}\|_{L^\infty} \le \eps/2.
\end{equation*}
We observe that for any $s\ge\eps$, this implies $\{ \hh > s \} \subset \{ \overline{h} > \eps/2 \}$, and so $\mathcal{H}^1(\widetilde I_s) \le \mathcal{H}^1(\{ \overline{ h} > \eps/2 \}) \le \mathcal{H}^1(\supp \overline{h}) = C$ for any $s \ge \eps$. From this follows that if $I_s$ denotes the largest connected component of $\{h_V>s\}$, where $h_V$ is obtained by the inverse rescaling of $\hh$, then
\begin{equation}\label{sizeIy}
\mathcal{H}^1( I_s)\le C V^{2/5} \qquad \qquad \forall s \ge \eps V^{3/5}.
\end{equation}

Now we claim that for some $t \in [2\eps V^{3/5},3\eps V^{3/5}]$, we have
\begin{equation}\label{decayener}
 \|u_V(\cdot, t)\|_{\dot H^{1/2}(I_t\times\{t\})}^2 :=  \min_{v(\cdot,t) = u_V(\cdot,t)} \int_{I_t\times [t,+\infty)} |\nabla v|^2 \dxy \les V^{4/5} \eps^{-1} \exp\lt(-C \eps V^{1/5}\rt).
\end{equation}
Indeed, fix $s \ge \eps V^{3/5}$. Then since $u_V$ is the minimizer of the Dirichlet energy, it satisfies the Laplace equation with Neumann boundary conditions at the upper part of the boundary. Denoting by $\bar u_V$ the average value of $u_V$ on $I_s \times \{s\}$, $\Omega_V^s:=\Omega_{h_V} \cap (I_s\times[s,+\infty))$, and using H\"older's and Poincar\'e's inequality, we get
\begin{align*}
 \int_{\Omega_V^s} |\nabla u_V|^2 \dxy &=\int_{I_s \times \{s\}} (u_V-\bar u_V) \partial_y u_V \dx\\
& \le \lt(\int_{I_s\times \{s\}} |u_V-\bar u_V|^2 \dx \rt)^{1/2} \lt( \int_{I_s\times \{s\}} \lt(\partial_y u_V\rt)^2 \dx \rt)^{1/2}\\
&\le \frac{\mathcal{H}^1(I_s)}{\pi} \lt( \int_{I_s\times \{s\}} \lt(\partial_x u_V\rt)^2\dx \rt)^{1/2} \lt( \int_{I_s\times \{s\}} \lt(\partial_y u_V\rt)^2 \dx \rt)^{1/2}\\
&\le \frac{\mathcal{H}^1(I_s)}{2\pi} \int_{I_s\times \{s\}} |\nabla u_V|^2 \dx \overset{\eqref{sizeIy}}{\le}C V^{2/5} \int_{I_s\times \{s\}} |\nabla u_V|^2 \dx.
\end{align*}
Since for $F(s):= \int_{\Omega_V^s} |\nabla u_V|^2 \dxy$, $s \ge \eps V^{3/5}$, this is equivalent to $F(s)\le -CV^{2/5} F'(s)$, an integration in $s$ yields that for $ s \ge \eps V^{3/5}$,
\[\int_{\Omega_V^s} |\nabla u_V|^2 \dxy\le \int_{\Omh} |\nabla u_V|^2  \dxy  \cdot \exp\lt(-C \frac{s-\eps V^{3/5}}{V^{2/5}}\rt)\le V \exp\lt(-C \frac{s-\eps V^{3/5}}{V^{2/5}}\rt).\]
In particular, the previous relation with $s = 2\eps V^{3/5}$ implies
\[ \int_{2\eps V^{3/5}}^{3\eps V^{3/5}} \int_{I_s} |\nabla u_V(x,s)|^2 \dx \ud s \le V \exp(-C \eps V^{1/5}), \]
and so there exists $t \in [2\eps V^{3/5},3\eps V^{3/5}]$ such that 
\[ \int_{I_t\times \{t\}} |\nabla u_V|^2 \dxy \le V^{2/5}\eps^{-1} \exp(-C \eps V^{1/5}).\]
Finally, we use Wirtinger's inequality 
\[ \| u_V - \bar u_V \|_{L^2(I_t\times\{t\})} \les \mathcal{H}^1(I_t) \cdot \| u_V \|_{\dot H^{1}(I_t\times\{t\})},\] where $\bar u_V$ denotes the average of $u_V$ on $(I_t\times\{t\})$, together with 
\[ \| u_V \|_{\dot H^{1/2}(I_t\times\{t\})}^2 \le \| u_V - \bar u_V\|_{L^2(I_t\times\{t\})} \cdot \| u_V \|_{\dot H^{1}(I_t\times\{t\})}\] 
(see, e.g., \cite[Eq. (9)]{Conti}), to get $\| u_V \|_{\dot H^{1/2}(I_t\times\{t\})} \les \mathcal{H}^1(I_t)^{1/2} \cdot \| u_V \|_{\dot H^{1}(I_t\times\{t\})}$ and thus \eqref{decayener}.

As a final step of the proof, for $s\ge 3\eps$ we want to construct $(\widetilde u, \widetilde h)$, a competitor for $(\widetilde u_V,\hh)$. Outside of $\widetilde I_s$ let $\widetilde h := \widetilde h_V$, and in $\widetilde I_s$ let $\widetilde h := \bar h_s$. Take then $\widetilde u$ to be equal to $u_V$ outside of $\Omega_V^t$ and to the restriction of a minimizer of
\[\min_{u(\cdot,t) =u_V(\cdot,t)} \int_{I_t\times[t,+\infty)} |\nabla u|^2 \dxy \]
elsewhere. 
By minimality of $(\widetilde h_V, u_V)$ we infer that
\[\int_{\widetilde I_s} \widetilde h_V'^2 \dx + V^{-4/5}\int_{\Omega_V^s} |\nabla u_V|^2 \dxy\le \int_{\widetilde I_s} \widetilde h'^2 \dx + V^{-4/5}\|u_V(\cdot, t)\|_{\dot H^{1/2}(I_t\times\{t\})}^2,
\]
hence
\[\int_{\widetilde I_s} \lt(\widetilde h_V'^2- \widetilde h'^2\rt) \dx \le V^{-4/5}\|u_V(\cdot, t)\|_{\dot H^{1/2}(I_t\times\{t\})}^2 \overset{\eqref{decayener}}{\les}  \eps^{-1} \exp\lt(-C \eps V^{1/5}\rt).
\]
Using Proposition \ref{prop:exp1} we conclude the proof.
\end{proof}

}

\section{The large slope approximation}\label{sec:largeslope}
If at large volume, a compact island forms, then $|h'|$ is expected to be large on its support. In this case, the small slope approximation $\sqrt{1+|h'|^2}-1\textcolor{red}{\sim} h'^2$ from Section \ref{sec:smallslope} might not be appropriate, and we rather work with the large slope approximation
\begin{eqnarray}
\sqrt{1+h'^2}-1\sim |h'|,
 \label{eq:largeslope}
\end{eqnarray}
i.e., we consider now the functional
\begin{equation}\label{eq:energylargeslope}
{{\int_{\Omega_h}|\nabla u|^2\dxy+\int_\R|h'|.}}
\end{equation}
Note that this approximation comes along with a loss of regularity of $\Omega_h$ since for low-energy configurations, $h$ is no longer bounded in $H^1(\R)$ but only in $BV(\R)$. Hence, we consider the relaxation of the energy as determined in the case of compact support in \cite{FFLM:07}. We follow the notation of \cite{FFLM:07, BC:02,GolZwick}. If $h:\R\to \left[0,+\infty\right)$ is lower semicontinuous {{(l.s.c.)}}, then we denote the pointwise variation of $h$ by
\[\int_\R|h'|:=\textup{Var}\ h:= \sup\left\{ \sum_{i=1}^n \left|h\left(x_i\right)-h\left(x_{i-1}\right)\right| : \quad x_1<\dots<x_n \right\}\ .\]
If $\textup{Var}\ h$ is finite then  $h$ is said to be of bounded pointwise variation (see \cite{AFP:00}). For a function $h$ of bounded  pointwise variation, set
\begin{align*}
 h^{-}\left(x\right)&:=\min \left\{h\left(x^+\right),h\left(x^-\right) \right\}=\liminf_{z\to x} h\left(z\right),\\
h^{+}\left(x\right)&:=\max \left\{h\left(x^+\right),h\left(x^-\right) \right\}=\limsup_{z\to x} h\left(z\right),
\end{align*}
where $h\left(x^{\pm}\right):= \lim_{z\to x^{\pm}} h\left(z\right)$.
We denote by $\Gamma_{cuts}$ the at most countable collection of vertical cuts,
\[\Gamma_{cuts}:=\left\{\left(x,y\right) \ : \ x\in S\left(h\right), \ h\left(x\right)\le y \le h^-\left(x\right)\right\} \]
 where $S\left(h\right):=\left\{x : h\left(x\right)<h^-\left(x\right)\right\}$. 
For $h$ lower semicontinuous and of bounded pointwise variation{{,}} we set
\[E_V(u,h):=\int_{\Omega_h}|\nabla u|^2\dxy\mbox{\quad and\quad}S_V^{{(\ell)}}(h):=\int_{\R}|h'|+2\mathcal{H}^1(\Gamma_{\text{cuts}}). \]
Note that $\mathcal{H}^1(\Gamma_{\text{cuts}})=0$ for locally Lipschitz functions $h$. If for a sequence of bounded energy we restrict ourselves to a compact set, we are in the situation of \cite{FFLM:07,GolZwick}, and we obtain a local compactness result by \cite[Lemma 2.1, Proposition  2.2, and Theorem 2.8]{FFLM:07}. The result becomes a global result, if we have strong $L^1$-convergence of $\{h_n\}$, which, in turn, follows from tightness. 
\begin{prop1}
 \label{prop:lsc}
Assume $(u_n,h_n)$ is admissible for \eqref{eq:energylargeslope} with $\int_\R h_n(x)\dx=V$. Then there exists a subsequence $(h_n, u_n)$ (not relabeled) such that $h_n\rightarrow h$ in $L^1_{\text{loc}}(\R)$ with 
\[h(x):=\inf\{\varliminf_{n\rightarrow\infty}h_n(x_n) \, : \, x_n\to x\}. \]
It holds that $\int_{\R}h\dx\leq V$. 
Further, $\R^2\setminus\Omega_{h_n}\rightarrow\R^2\setminus\Omega_h$ in the local Hausdorff topology, and  $u_n\rightharpoonup u$ in $H^1_{\text{loc}}(\Omega_h)$. If $\{h_n\}$ is tight, then the convergences hold globally, $\int_\R h\dx=V$, and
\[\varliminf_{n\rightarrow\infty} S_V^{{(\ell)}}(h_n)+E_V(u_n,h_n)\geq S_V^{{(\ell)}}(h)+E_V(u,h). \]
Conversely, if $h$ is a lower semicontinuous, nonnegative function of bounded pointwise variation, with $\int_\R h\dx=V$, then there exists a sequence of compactly supported Lipschitz functions $h_n$ with $\int_\R h_n(x)\dx=V$, and $u_n\in H^{1}(\Omega_{h_n})$ such that $h_n\rightarrow h$ in $L^1(\R)$, $u_n\rightharpoonup u$ in $H^1_{\text{loc}}(\Omega_h)$, and 
\[S_V^{{(\ell)}}(h)+E_V(u,h)\geq\varlimsup_{n\rightarrow\infty}S_V^{{(\ell)}}(h_n)+E_V(u_n,h_n). \]
\end{prop1}
{{For $V>0$, we set (see \cite[(2.8)]{FFLM:07})
\begin{eqnarray*}
 {\mathcal{A}}^{\ell}_V:=\left\{ (u,h):\ h\in BV(\R),\ h \text{\ l.s.c},\  h\geq 0, \ u\in H^1(\Omega_h),\ \int_\R h\, \dx=V,\,  u(x,0)= x \ \text{if\ }x\in\supp h\right\},
\end{eqnarray*}
and 
\begin{eqnarray}\label{eq:Fl}
\F (V):=\inf\left\{\int_{\Omega_h}|\nabla u|^2\dxy+\int_\R|h'|:\,
(h,u)\in {\mathcal{A}}^{\ell}_V\right\}.
\end{eqnarray}
}}
\subsection{Basic properties}
Note that the rescaling property for the small slope approximation from Lemma \ref{rescale} carries over to the case of the large slope approximation \eqref{eq:largeslope}. Precisely, for a lower semicontinuous function $h\in {{BV}}(\R)$, $u\in H^1(\Omega_h)$ and $\lambda>0$, consider the rescaled quantities $h_\lambda\in {{BV}}(\R)$ and $u_\lambda\in H^1(\Omega_h)$ given by $h_\lambda(x):=\frac{1}{\lambda}h(\lambda x)$ and $u_\lambda(x,y):=\frac{1}{\lambda}u(\lambda x,\lambda y)$. Then $\int_{\R}h\dx=\lambda^2\int_{\R}h_\lambda(x)\dx$, and
\[\int_{\Omega_h}|\nabla u|^2\dxy+\int_{\R}|h'|+2\mathcal{H}^1(\Gamma_{\text{cuts}})=\lambda^2\int_{\Omega_{h_\lambda}}|\nabla u_\lambda|^2\dxy+\lambda\int_\R|h_\lambda'|+2\lambda\mathcal{H}^1(\Gamma_{\text{cuts}}). \]
Further, 
we again have $F_\ell(V)\leq V$ by considering a sequence $h_n(x)=\frac{V}{n}\chi_{(0,n)}$ and $u_n(x,y)=x$ in $\Omega_{h_n}$. 
Hence, some results from the previous sections carry over essentially verbatim to the large slope approximation. We collect some properties in the following proposition.
\begin{prop1}
\label{rem:basicsls}
\begin{itemize}
\item[(i)] $\F$ is concave, and hence in particular locally Lipschitz continuous. At every point of differentiability,
\[F_\ell'(V)=\frac{1}{V}(E_V+\frac{1}{2}S_V). \]
Minimizers for a fixed volume have the same surface and the same elastic energy. Further,
\[
 \overline{\lim}_{\eps\rightarrow 0}\frac{F_\ell(V+\eps)+F_\ell(V-\eps)-2F_\ell(V)}{\eps^2}\leq-\frac{S_V}{4V^2}. 
\]
\item[(ii)] For every $V>0$, {{if a minimizer exists, i.e., if there is $(u^\ast,h^\ast)\in{\mathcal{A}}^{\ell}_V$ such that
$\F(V)=E_V(u^\ast,h^\ast)+S_V^{(\ell)}(h^\ast)$}}, then $\overline{\{h{{^\ast}}>0\}}$ is connected. 


\item[(iii)] For every $V>0$,
\begin{multline*}\F(V)=\inf\left\{\int_{\Omega_h}|\nabla u|^2 \dxy+\int_{\R}|h'|\ :\ {{(u,h)\in\mathcal{A}_V^{\ell},\ \partial_yu\equiv 0\text{\ or\ }}}\int_{\Omega_h}(\partial_y u)^2 \dxy=\frac{1}{4}\int_{\R}|h'|\right\}
\end{multline*}
\item[(iv)] We have 
\[\|h\|_{L^\infty(\R)}\leq\int_{\R}|h'|. \]
\end{itemize}
\end{prop1}
\begin{proof}
 \begin{itemize}
\item[(i)] See proofs of Propositions \ref{concavss} and \ref{prop:diff}.
\item[(ii)] See proof of Proposition \ref{connect}.
  \item[(iii)] Consider the anisotropic volume-preserving scaling from the proof of Lemma \ref{internvar}, i.e., $h_\lambda(x):=\lambda h(\lambda x)$, and $u_\lambda(x):=\frac{1}{\lambda}u(\lambda x,\frac{1}{\lambda}y)$. Then
\begin{eqnarray*}
 \int_{\Omega_h}|\nabla u_\lambda|^2 \dxy+\int_{\R}|h_\lambda'|=\int_{\Omega_h}\left((\partial_x u)^2+\frac{1}{\lambda^4}(\partial_y u)^2\right)\dxy+\lambda\int_{\R}|h'|,
\end{eqnarray*}
{{
and for $\partial_yu\not\equiv 0$, minimization in $\lambda$ yields the claim.
}}
  \end{itemize}
\end{proof}

\subsection{Scaling law
}
The scaling law for $\F$ can be derived arguing along the lines of the proof of Proposition \ref{prop:scaling}.
\begin{prop1}
\label{prop:scalinglargeslop}
 There is a constant $c_0>0$ such that for all $V>0$, we have
\[c_0\min\{V,\,V^{2/3}\}\leq\F(V)\leq\frac{1}{c_0}\min\{V,\,V^{2/3}\}
      .
 \]
Further, there is a constant $c>0$ with the following property: If $V$ is large enough, and $(u,h)$ is admissible for \eqref{eq:energylargeslope} with $F(u,h)\leq\frac{1}{c_0}\min\{V,\,V^{2/3}\}$, then $\sup h\geq cV^{2/3}$.
\end{prop1}
\begin{proof}
 We prove the upper bound first. If $V \le 1$, for ${{L}}
>0$ we set $h:=V{{L}}
^{-1}\chi_{[0,{{L}}
]}$ and $u(x,y):=x$ for all $(x,y)\in\Omega_h$. Then 
\begin{eqnarray*}
 \F(V)\leq \int_0^{{L}}
\int_0^{V{{L}}
^{-1}}|\nabla u|^2\dxy+\int_{\R}|h'|= V + 2V{{L}}
^{-1}. 
\end{eqnarray*}
Since ${{L}}
 > 0$ can be chosen arbitrary large, we have that $\F(V) \le V$. 

In the case $V\ge 1$, choose ${{L}}
:=V^{1/3}$, $h:=\frac{V}{{{L}}
}\chi_{[0,{{L}}
]}$, and let $u$ be a minimizer of the Dirichlet energy on $\Omega_h=[0,{{L}}
]\times[0,h]$ subject to $u(x,0)=x$. 
Since ${{L}}
= V^{1/3}\ll V^{2/3}\sim \frac{V}{{{L}}
}$, we have \[\int_{\Omega_h}|\nabla u|^2 \dxy\sim{{L}}
^2\sim V^{2/3}\mbox{\quad and\quad}\int_{\R}|h'|= 2V^{2/3}. \]
The lower bound together with the estimate on $\sup h$ can be obtained by repeating the proof of Proposition \ref{prop:scaling} replacing the estimate $S_1\ges\frac{V_1^2}{t_1^3}$ by $S_1\ges \frac{V_1}{t_1^2}$. .
\end{proof}
\begin{remark1}
 Just as for $F_s$, using (iv) of Proposition \ref{rem:basicsls}, we can obtain that for an almost minimizer, $\sup h\sim V^{2/3}$, and the size of the support of $h$ is at least $\sim V^{1/3}$.
\end{remark1}

\begin{remark1}
 If we rescale $V:=e_0^4d$, and $\F(e_0,d)=\frac{1}{e_0}\F(V)$, we find
\[\min\F(e_0,d)\sim \min\{e_0^2d,\,e_0^{2/3}d^{2/3}\}. \]
\end{remark1}
Over the last years, much work has been devoted to the analysis of island formation in a compact setting (see \cite{FFLM:07,FM:12,GolZwick}). Precisely, assuming that $h:[0,1]\rightarrow[0,\infty)$ is Lipschitz, set
\begin{multline}\label{eq:golzwi}
\tilde{F}(d):=\inf\left\{\int_{\Omega_h}|\nabla u|^2\dxy+\int_0^1 (\sqrt{1+(h')^2}-1)\dx\ :\right.\\
\left. {{h\in W^{1,\infty}(\R), \,u\in H^1(\Omega_h)}},\,\int_0^1h\dx=d,\,u(x,0)=e_0{{x}}, h(0)=h(1)=0\right\}. 
\end{multline}
In this case, the surface energy is always bounded below by $\min\{d,\,d^2\}$ since by the compact support there is a point $x^\ast\in(0,1)$ with $h(x^\ast)\geq d$, and thus 
\begin{eqnarray}
&& \int_0^1 (\sqrt{1+h'^2}-1)\dx=\int_0^{x^\ast} (\sqrt{1+h'^2}-1)\dx+\int_{x^\ast}^1 (\sqrt{1+h'^2}-1)\dx\nonumber\\
&\geq&\sqrt{1+(2d)^2}-1\geq\begin{cases}
                                                          cd&\mbox{if $d$ is large}\\                                                                                      
                                                                                     cd^2&\mbox{if $d$ is small.}                                                          \end{cases}
\label{eq:compsurf}
\end{eqnarray}
Consequently, for compact support, the scaling law is the following. 
\begin{prop1}\label{rem:scalingcompls}
The following holds
\begin{eqnarray}
 \tilde{F}(d)\sim \max\{\min\{d,\,d^2\},\,\min\{de_0^2,\,d^{2/3}e_0^{2/3}\}\}.
\label{eq:scalingcomp}
\end{eqnarray}
\end{prop1}
\begin{proof}
 We prove the upper bound first.
\begin{itemize}
 \item[(i)] If $d\lesssim 1$, let $\Omega_h$ be a triangle of length $1$ and height $2d$, and set $u(x,y)=e_0x$ in $\Omega_h$. Then, since $|h'|$ is small,
\[\tilde{F}(u,h)\lesssim e_0^2d+d^2\lesssim\begin{cases}
                                            d^2&\mbox{if\ }e_0^2\lesssim d\\
e_0^2d&\mbox{if\ }e_0^2\gtrsim d.
                                           \end{cases}
\]
\item[(ii)] If $1\lesssim e_0^2\lesssim d$, let $\Omega_h$ be close to a rectangle of length $1$ and height $d$, i.e., $\supp h\subset (0,\eps)\cup(1-\eps,1)$, and let $u$ be a minimizer of the Dirichlet energy in $\Omega_h$ subject to the boundary condition $u(x,0)=e_0x$. Then, since $|h'|$ is large on its support,
\[\tilde{F}(u,h)\lesssim e_0^2+d\lesssim d. \] 
\item[(iii)] If $\frac{1}{e_0^4}\lesssim d\lesssim e_0^2$, let $\Omega_h$ be a Lipschitz approximation of a rectangle of length ${{L}}
\sim\frac{d^{1/3}}{e_0^{2/3}}\lesssim 1$, and height $h\sim d^{2/3}e_0^{2/3}$, and let $u$ be a minimizer of the Dirichlet integral in $\Omega_h$ subject to the boundary condition $u(x,0)=e_0x$. Then, since $|h'|$ is large on its support,
\[\tilde{F}(u,h)\lesssim e_0^2{{L}}
^2+d^{2/3}e_0^{2/3}\sim e_0^{2/3}d^{2/3}. \]
\end{itemize}
The lower bound follows from \eqref{eq:compsurf} and a proof similar to that of Proposition \ref{prop:scaling} since in case $V_1>2 t_1^2$, the slope used to estimate the surface energy in large.
\end{proof}
 The scaling law \eqref{eq:scalingcomp} resembles essentially results from \cite{GolZwick}, where a model without normalization of the surface energy has been considered, i.e.,
\begin{multline*}
\tilde{F}_2(d):=\inf\left\{\int_{\Omega_h}|\nabla u|^2\dxy +\int_0^1 \sqrt{1+(h')^2}\dx :\right.\\
\left. {{h\in W^{1,\infty}(\R), \,u\in H^1(\Omega_h)}},\,\int_0^1h\dx=d,\,u(x,0)=e_0{{x}}, h(0)=h(1)=0\right\}. 
\end{multline*}
 There, the scaling law turns out to be 
\[\inf\tilde{F}_2\sim\max\{1,\,d, \, e_0^{2/3}d^{2/3}\}. \]
We note that the proof uses essentially the large slope approximation $\sqrt{1+h'^2}\geq\max\{1,\,|h'|\}$, and small slopes (which are likely only for small volumes) are not seen by the first order behavior. 

\subsection{Existence  of minimizers}
In case of the small slope approximation, non-existence of minimizers follows from the fact that there is a regime of volumes for which $F_s(V)=V$ (see Proposition \ref{prop:nonexist}). In case of the large slope approximation, we can prove only a weaker statement in this direction.
\begin{prop1}\label{prop:smallvollargeslope}
 For every $\delta>0$ there is $\widetilde{V}(\delta)>0$ such that for every $0 < V < \widetilde{V}$,  
\[(1-\delta)V\le F_\ell(V) \le V.\]
In particular,
\[\lim_{V\rightarrow 0}\frac{1}{V}F_\ell(V)=1. \]
\end{prop1}
\begin{proof}
First, the upper bound $F_\ell(V) \le V$ is showed in the proof of Proposition~\ref{prop:scalinglargeslop}. To prove the 
{{other estimate, }}
we refine the argument 
{{from}}
the 
{{proof of the}} 
lower bound for the minimal energy {{in Proposition~\ref{prop:scalinglargeslop}}}. 

Let $h$ be any locally Lipschitz function with $\int_{\R}h\dx=V$ and $0 < \delta < 1$ be fixed. 
Let $\lambda>0$ be a parameter which will be chosen below. Pick $x_1\in\R$ and define $t_1:=\max\{t>0:\,[x_1,x_1+\lambda t]\times[0,t]\subset\Omega_h\}$. Set $V_1:=\int_{x_1}^{x_1+t_1}h\dx > 0$.

Let us first prove that there exists ${{\widetilde{V}:=}}\widetilde{V}(\delta)$ such that if $V_1 <  \widetilde{V}$, then 
\begin{equation}\label{limitlb}
 \int_{x_1}^{x_1+t_1} \int_0^{h(x)} |\nabla u|^2 \dx \dy + \int_{x_1}^{x_1+t_1} |h'| \ge (1-\delta)V_1.
\end{equation}
{{Assume for the sake of contradiction that}}
\begin{eqnarray}\label{eq:48AFSOC}
\int_{x_1}^{x_1+t_1} \int_0^{h(x)} |\nabla u|^2 \dx \dy + \int_{x_1}^{x_1+t_1} |h'| < (1-\delta)V_1.
\end{eqnarray}
 Since $h$ is locally Lipschitz,
  we have that $\max_{[x_1,x_1+t_1]} h - \min_{[x_1,x_1+t_1]} h \le \int_{x_1}^{x_1+t_1} |h'|$. Then $\max_{[x_1,x_1+t_1]} h \ge V_1/t_1$ and $\min_{[x_1,x_1+t_1]} h = \lambda t_1${{, which implies that} }
  $V_1/t_1 - \lambda t_1 < (1-\delta)V_1$, and subsequently 
  \begin{eqnarray}\label{eq:48est1}
  V_1(1 - (1-\delta)t_1) < \lambda t_1^2. 
\end{eqnarray}
Next, by Lemma~\ref{LeDretRaoult}{{,}} there exists $\psi = \psi(\lambda)$ such that 
\begin{equation}\label{limitlb2}
 \int_{x_1}^{x_1+t_1} \int_0^{\lambda t_1} |\nabla u|^2 \dxy \ge \psi \lambda t_1^2,\text{\quad and \ $\psi \to 1$ as $\lambda \to 0$. }
\end{equation}
In particular, we can choose $\lambda=\lambda(\delta) > 0$ small enough such that $\psi(\lambda) > 1-\delta$. 
{{By~}}\eqref{limitlb2}, {{the assumption~\eqref{eq:48AFSOC}}} implies $\psi \lambda t_1^2 < (1-\delta)V_1$. We combine this with{{~\eqref{eq:48est1}}}
 to obtain
\begin{equation}\label{limitlb3}
 t_1 > \frac{1-\frac{1-\delta}{\psi}}{1-\delta} > 0. 
\end{equation}
Now we define $\widetilde{V}(\delta) := \lambda
 \left(\frac{1-\frac{1-\delta}{\psi}}{1-\delta}\right)^2 > 0$ and observe that
 ~\eqref{limitlb3} implies $V_1 \ge \lambda t_1^2 \ge \widetilde{V}$, a contradiction 
 {{to}} the assumption $V_1 < \widetilde{V}$.

If we choose $\widetilde{V}$ as above, we can continue the same way as in the proof of Proposition~\ref{prop:scaling} 
and use that $V_i \le V < \widetilde{V}$ to get
\begin{equation}
 F_\ell(V) \ge (1-\delta) \sum_{i} V_i = (1-\delta)V. 
\end{equation}
\end{proof}
 
This slightly weaker statement still allows us to derive an analogue to Lemma \ref{decrbeta}. We define $\beta_\ell$ by $F_\ell(V)=:\beta_\ell(V)V$. Note that again $\beta_\ell(V)<1$ for large $V$. 
\begin{lemma1}
 \label{lem:decrbetal}
The function $\beta_\ell$ is 
strictly
decreasing in $\{F_\ell(V)<V\}=\{\beta_\ell<1\}$.
\end{lemma1}
\begin{proof}
 Assume for a contradiction that there are $V_0<V_1$ such that $F_\ell(V)=\beta_0V<V$ for all $V\in[V_0,V_1]$ with some $0<\beta_0<1$. By Proposition \ref{prop:smallvollargeslope}, there exists $\tilde{V}>0$ such that $F_\ell(V)\geq\frac{1}{2}(1+\beta_0)V$ for all $V<\tilde{V}$. As in the proof of Lemma \ref{decrbeta}, concavity of $F_\ell$ and $F_\ell(0)=0$ imply that $F_\ell(V)=\beta_0V$ for $V\in(0,\tilde{V})$, which yields a contradiction.
\end{proof}
Proceeding along the lines of Section \ref{sec:smallslope}, we show tightness of minimizing sequences for $V\in\{\beta_\ell<1\}$.
\begin{lemma1}
 \label{lem:cuttingl}
Let $V\in\{\beta_\ell<1\}$ and $\delta>0$. Then there exist $\ell=\ell(V,\delta)>0$ and $C(V,\delta)>0$ with the following property: If $(u,h){{\in\mathcal{A}_V^{\ell}}}$ 
with 
$\eps:=E_V(u,h)+S_V(h)-F_\ell(V)\leq C(V,\delta)$, and $x_0<x_1$ with $x_1-x_0=\ell$, then 
\[\int_{-\infty}^{x_0}h\dx\leq\delta\mbox{\qquad or\qquad}\int_{x_1}^\infty h\dx\leq\delta. \]
\end{lemma1}
\begin{proof}
 We proceed along the lines of the proof of Lemma \ref{lm:cutting}. {{Fix $V\in\{\beta_\ell<1\}$ and $\delta>0$. Let 
 \[ 0<C(V,\delta)\leq\frac{\delta}{2}\left(\beta_\ell\lt(\frac{V}{2}\rt)-\beta_\ell(V)\right),\text{\qquad and\quad}0<\alpha<\frac{\delta}{36}\left(\beta_\ell\lt(\frac{V}{2}\rt)-\beta_\ell(V)\right)\] 
 }}
 be 
 such that $V/(3\alpha)=n\in\N$. Define $\ell:=\alpha^{-1}V$, and consider an interval $[x_0,x_1]$ of length $\ell$, and write it as a disjoint union of $3n$ intervals of length $1$. As in the proof of Lemma \ref{lm:cutting}, taking into account the different surface energy term, there is an interval $I$ such that
\[\int_Ih\,\dx{{\leq 3\alpha}}\mbox{\qquad and \qquad}\int_I|h'|\leq3\alpha, \]
which implies 
\[3\alpha\geq\int_I|h'|\geq\sup_Ih-\inf_Ih. \]
Since also $3\alpha\geq\int_Ih\dx\geq\inf_Ih$, we have
\[\sup_Ih=(\sup_I h-\inf_Ih)+\inf_Ih\leq3\alpha+3\alpha{{=6\alpha}}. \]
{{Following the lines of the proof of Lemma~\ref{lm:cutting}, w}}e make two cuts in $I$ such that the profile is separated into two pieces of volumes $V_0<V_1$, with $V=V_0+V_1+V_{lo}$, with  ${{V_{lo}\leq
16\alpha}}$
such that the surface energy is increased by at most $S_{cut}\leq 2\sup_Ih\leq {{12\alpha}}$.
Then as in \eqref{v0}, $V_0\left(\beta_\ell(V/2)-\beta(V)\right)\leq\eps+V_{lo}+S_{cut}$, and thus, 
\[V_0\leq\frac{\eps+V_{lo}+S_{cut}}{\beta_\ell(V/2)-\beta_\ell(V)}\leq {{\frac{\eps+6\alpha+12\alpha}{\beta_\ell(V/2)-\beta_\ell(V)}<\delta. }}
\]
The proof is concluded as in Lemma \ref{lm:cutting}.
\end{proof}
As worked out in Proposition \ref{existminss}, Lemma \ref{lem:cuttingl} implies tightness (up to translations) of minimizing sequences. Using lower semicontinuity of the energy we obtain existence of minimizers:
\begin{prop1}
 \label{prop:existminls}
Let $V$ be such that $F_\ell(V)<V$. Then there exists a minimizer (u,h) of {{\eqref{eq:Fl}}}.
\end{prop1}

\subsection{Asymptotic behavior}
We proceed along the lines of Section \ref{sec:asympss}. According to the slightly different scaling law (see Proposition  \ref{prop:scalinglargeslop}), we define for admissible $h$ with $\int_\R h\,\dx=V$ the rescaled quantities by $\tilde{h}(x):=V^{-2/3}h(V^{1/3}x)$, and 
\[G_V(\tilde{h}):=V^{-2/3}\left(E_V(u,h)+S_V(h)\right)=V^{-2/3}\int_{\Omega_h}|\nabla u|^2\dxy+\int_{\R}|\tilde{h}'|. \]
Note that $\int_{\R}\tilde{h} \dx=1$.
\begin{teorema1}
 \label{th:gammals}
 For every sequence $V_n\rightarrow\infty$, let $(u_{V_n},h_{V_n})$ be a sequence of minimizers {{of \eqref{eq:Fl}}} . Then there exists a subsequence such that the rescaled profile functions $\tilde{h}_{V_n}$ converge to $h$ in $L^1(\R)$, which minimizes
\begin{equation}\label{Gls}G(h):=\left(\inf_{u(x,0)=x}\int_{\{h>0\}\times[0,+\infty)} |\nabla u|^2 \dxy\right)+\int_\R |h'|+2\mathcal{H}^1(\Gamma_{\text{cut}})
\end{equation}
 subject to the constraint $\int_\R h \dx=1$. Moreover, $\R^2\backslash \Omega_{\tilde h_{V_n}}$ converges in Hausdorff topology to $\R^2\backslash \Omega_{ h}$.
\end{teorema1}
\begin{proof}
 We first observe that $G_{V_n}(\tilde{h}_{V_n})\leq C$, and that (up to translations) the sequence $\tilde{h}_{V_n}$ is tight (see proofs of Theorem \ref{th:asympss} and Lemma \ref{lem:cuttingl}). Thus, a subsequence of $\tilde{h}_{V_n}$ converges to $h$ in $L^1(\R)$. Thanks to the bound on the surface energy, we also have local Hausdorff convergence of $\R^2\backslash \Omega_{\tilde h_{V_n}}$ to $\R^2\backslash \Omega_{ h}$ which then improves to Hausdorff convergence thanks to Proposition \ref{rem:basicsls}  (iv).
  By the lower semicontinuity, we obtain the lower bound for the surface energy, i.e.,
 \[\varliminf_{n\rightarrow\infty}\int_{\R}|\tilde{h}_{V_n}|\geq \int_{\R}|h'|+\mathcal{H}^1(\Gamma_{cut}). \] 
For the elastic energy, we proceed as in \cite[Proposition 4.3]{GolZwick}. By the Hausdorff convergence of $\Omega_{\tilde{h}_V}$ to $\Omega_h$, we have convergence of the ``boundary layers'' $\{\tilde{h}_{V_n}>0 \text{\ and\ }h>0\}\times [0,\infty)$ to $\{h>0\}\times [0,\infty)$ in the local Hausdorff topology. Changing variables, $x=V_n^{1/3}\hat{x}$, $y=V_n^{1/3}\hat{y}$, $u_{V_n}=V_n^{-1/3}\hat{u}_{V_n}$, there is a subsequence such that $\hat{u}_{V_n}\rightharpoonup\hat{u}$ locally weakly in $H^1(\{h>0\}\times[0,\infty))$. Thus as in \cite[Proposition 4.3]{GolZwick},
\[V_n^{-2/3}\int_{\Omega_{h_{V_n}}}|\nabla u_{V_n}|^2\dxy\geq  \int_{\Omega_{\tilde{h}_{V_n}}\cap[\{h>0\}\times[0,\infty)]}|\nabla \hat{u}_{V_n}|^2\dxy\geq \int_{\{h>0\}\times[0,\infty)}|\nabla\hat{u}_{V_n}|^2\dxy.\]
The lower bound follows. {{For the upper bound construction, it is enough to consider $\tilde{h}_{V_n}=h$.}} We conclude by $\Gamma$-convergence, as in the proof of Theorem \ref{th:asympss}.
\end{proof}

\begin{remark1}
\label{rem:limitenergyls}
 We note that
\[G(h):= C_W \sum_{i\in \N} (b_i-a_i)^2+\int_\R |h'|+2\mathcal{H}^1(\Gamma_{\text{cuts}}),\]
  where the  intervals $(a_i,b_i)$ are the connected components of $\{h>0\}$, 
and $C_W$ is as in \eqref{eq:CW}.
\end{remark1}
We next study the minimizer of the limiting functional.
\begin{prop1}
\label{prop:minls}
 The minimization problem \eqref{Gls} admits a unique minimizer up to translations, namely the rectangle with base $2^{1/3}C_W^{-1/3}$.
\end{prop1}
\begin{proof}
 The proof follows as in \cite[Proposition 4.5]{GolZwick}, and we briefly sketch it only for the reader's convenience. The optimal function $h$ is of the form $h=\sum_{i=1}^Nh_i\chi_{(a_i,b_i)}$ with $\sum_{i=1}^Nh_i(b_i-a_i)=1$, and $a_i<b_i<a_{i+1}$ since the rectangle minimizes $\int|h'|$ among profiles with given volume. We set $\ell_i:=b_i-a_i$. Then by Remark \ref{rem:limitenergyls}, the energy is given by
\[C_W\sum_{i=1}^N \ell_i^2+2\sum_{i=1}^Nh_i.\]
Assume that there are two connected components, say, of lengths $\ell_1\geq\ell_2>0$, then for $\eta\in[-h_1,\frac{\ell_2}{\ell_1}h_2]$, we consider the volume-preserving variation of $h$ changing $h_1$ and $h_2$ to $h_1+\eta$ and $h_2-\eta\frac{\ell_1}{\ell_2}$, respectively. The minimality condition then gives $\ell_1=\ell_2$, from which we deduce that $\ell_i=\ell\equiv const$ for every $i=1,\dots,N$. The minimization problem then reduces to minimizing $C_WN\ell^2+\frac{2}{\ell}$ subject to $N\in\N$, which yields $N=1$ and $\ell=2^{1/3}C_W^{-1/3}$.
\end{proof}
To prove the exponential convergence, we have an analogue to Proposition \ref{prop:exp1}.
\begin{lemma1}
\label{exp1ls}
{{Let $L>0$ and $V>0$. C}}onsider $h_{\text{min}}:=\frac{V}{L}\chi_{[0,L]}$, which is the minimizer of $\int_{\R}|h'|$ 
{{in
\[X:=\{h\in BV(\R):  \ h\geq 0,\ \int_{\R} h(x)\dx=\int_0^Lh(x)\dx=V\}.\]}}
Then for all $h{{\in X}}$ 
we have
\[\int_\R|h'|-\int_\R|h_{\text{min}}'|\geq\frac{1}{L}\int_\R|h-h_{\text{min}}| \dx. \] 
\end{lemma1}
\begin{proof}
By rescaling {{the dependent and independent variables (see the proof of Lemma \ref{prop:interpol})}}, it suffices to consider $L=1$ and $V=1$, i.e., $h_{\text{min}}=\chi_{[0,1]}$ and $\int_{\R}|h_{\text{min}}'|=2$. By density, we may assume that $h\in W^{1,\infty}(\R)$. 
{{Then }}the function $h(x) - 1$ attain{{s}} its non-negative maximum at {{some}} $\overline{x}\in[0,1]$.  
We have $\int_0^1|h'| \geq 2h(\overline{x})$, and thus,
\begin{eqnarray*}
\int_0^1|h'|-2\geq 2(h(\overline{x})-1) \ge 2\int_0^1 (h-1)_{+}\dx=\int_0^1|h-h_{\text{min}}|\dx,
\end{eqnarray*}
where we used $\int_0^1 h \dx = 1$ to show that $2\int_0^1 (h-1)_+ \dx = \int_0^1 |h-1| \dx$. 
\end{proof}
Finally, we prove exponential convergence of a sequence of minimizing profiles of $G_V$ to the rectangle, which minimizes the limit functional (see Proposition \ref{prop:minls}). We denote by $I_s$ the largest connected component of $\{h_V>s\}$, and by $\tilde{I}_s$ the rescaled one. Given a function $\tilde{h}$, we denote by $\bar{h}_s$ the function that agrees with $\tilde{h}$ outside $I_s$, has the same volume as $\tilde{h}$, and is such that the surface energy term is minimized.
\begin{prop1}
 For every $\eps>0$ there exist constants $C_0$ and $C_1$ such that for every $V>\overline{V}$, and for every minimizer $\tilde{h}_V$ of $G_V$, for all $s\geq\eps$,
\[\|\tilde{h}_V-\overline{h}_s\|_{L^1(\tilde{I}_s)}\leq C_0\exp(-C_1V^{1/3}).\]
\end{prop1}
\begin{proof}
Let $\eps>0$. Since  (after possible translations), $\R^2\backslash \Omega_{\tilde h_{V}}$ converges in Hausdorff topology to $\R^2\backslash \Omega_{ h}$ and since $\overline{h}$ is a characteristic function,
\[\mathcal{H}^1(\tilde{h}>s)\leq C\qquad \forall s\geq\eps. \]
Hence, by rescaling, we obtain that for $s\geq \eps V^{2/3}$, the largest connected component $I_s$ of $\{h_V>s\}$ satisfies
\begin{eqnarray}
\label{eq:sizeIl}
\mathcal{H}^1(I_s)\leq CV^{1/3}.
\end{eqnarray}
By density, possibly slightly changing $\eps$, we may  instead of minimizers consider Lipschitz functions $\tilde{H}_V$ with $\|\tilde{H}_V-\tilde{h}_V\|_{L^1}\leq \exp(-C_1V^{1/3})$ and $E_V(u_V,H_V)+S_V(H_V)\leq E_V(u,h_V)+S_V(h_V)+\exp(-C_1V^{1/3})$.  
Similarly to the derivation of \eqref{decayener} we then obtain that there is some $t\in [2\eps V^{2/3},3\eps V^{2/3}]$ with
\begin{eqnarray}
\|u_V\|^2_{\dot{H}^{1/2}(I_t\times\{t\})}\leq V^{2/3}\eps^{-1}\exp(-C\eps V^{1/3}).
\label{eq:decayls}
\end{eqnarray}
Indeed, since $\tilde{H}_V$ is Lipschitz, and $u_V$ is a minimizer of the Dirichlet energy, we have as in the proof of Proposition \ref{expss},
\[\int_{\Omega_V^s}|\nabla u_V|^2\dxy\leq C\mathcal{H}^1(I_s)\int_{I_s\times\{s\}}|\nabla u_V|^2\dxy\leq CV^{1/3}\int_{I_s\times\{s\}}|\nabla u_V|^2\dxy, \]
i.e., $f(s):=\int_{\Omega_V^s}|\nabla u_V|^2$ for $s\geq\eps V^{2/3}$ satisfies $f(s)\leq -CV^{1/3}f'(s)$. Thus, for $s\geq\eps V^{2/3}$, we have $f(s)\leq V\exp\left(-C\frac{s-\eps V^{2/3}}{V^{1/3}}\right)$. In particular,
\[\int_{2\eps V^{2/3}}^{3\eps V^{2/3}}|\nabla u_V|^2\dxy\leq V\exp(-C\eps V^{1/3}), \]
and it follows by Wirtinger's inequality that there is some $t\in[2\eps V^{2/3},3\eps V^{2/3}]$ such that \eqref{eq:decayls} holds.
Now we choose as a competitor the function $\overline{h}_t$, which by definition agrees with $\tilde{H}$ outside of  $I_t$, and the corresponding optimal deformation $u$ with boundary data $u(x\cdot,t)=u_V(x,t)$ for all $x\in I_t$. Using the almost optimality, we conclude as in Proposition \ref{expss}, using Lemma \ref{exp1ls}. Note that the factor $V^{2/3}$ cancels with the rescaling factor of the elastic energy.
\end{proof}

\section{The scaling law in three space dimensions}\label{sec:3D}
{
\renewcommand{\F}{F^{3D}_s}
In this section we will identify the scaling law for the energy in the {{($2+1$)}} -dimensional setting. More precisely, we consider a $3D$ analog of~\eqref{funcss}: {{G}}iven $V \in (0,\infty)$, we define
\begin{eqnarray}\label{funcss3D}
 \F(V):=\inf
  \left\{\int_{\Omega_h} |\nabla u|^2 \dxyz + \int_{\R^2} |\nabla h|^2 \dx \dy : {{\ h\in W^{1,\infty}(\R),\ u\in H^1(\Omega_h)}},
 \right.\nonumber\\
 \left.\int_{\R^2} h \dx\dy=V,  \, u(x,y,0)=(x,y) {{\ \text{in }\supp h}}\right\},
\end{eqnarray}
where 
$\Omega_h := \{ \x:=(x,y,z) \in \R^3 : 0 < z < h(x,y) \}$.
\begin{Theorem}\label{th:scaling3D}
 There exists a positive constant $c$ such that for every $V > 0$
 \begin{equation*}
  c \min \{V,V^{6/7}\} \le \F(V) \le c^{-1} \min\{V,V^{6/7}\}.
 \end{equation*} 
\end{Theorem}

\begin{proof}
 To prove the upper bound, we consider two different constructions -- a thin layer and a pyramid.

 {\noindent \bf Thin layer construction:} For $\eps > 0$, let  $L$ be such that $\int_{\R^2} h \dx\dy = V$ for $h$ defined by $h(x,y) := \min{{\{}}\eps,\ [(L+\eps)-\max(|x|,|y|)]_+{{\}}}$, where by $f_+ := \max\{ 0, f \}$ denotes the positive part of $f$. 
 By setting $u(x,y,z) := (x,y)$ we get $\int_{\Omega_h} |\nabla u|^2 \dxyz = 2V$. We observe that $|\nabla h| = \chi_M$ with $M$ the difference of two concentric squares with sidelengths $2(L+\eps)$ and $2L$, and so
\begin{equation*}
  \int_{\R^2} |\nabla h|^2 \dx\dy = |M| = 4\eps (2L + \eps) \le 4(\sqrt{V\eps} + \eps^2),
 \end{equation*}
 where we used that $(2L)^2 \eps \le V$. Hence $F(V) \le 2V + 4(\sqrt{V\eps} + \eps^2)$ for arbitrarily small $\eps > 0$, which implies $F(V) \le 2V$. 
 
 {\noindent \bf Pyramid construction:} Let $L:=V^{2/7}$ and $H:=3V^{3/7}/4$. We define  \\
$h(x,y) := H \left( 1 - \max(|x|,|y|)/L\right)_+$ and $u(x,y,z) := (x,y) (1-z/L)_+$, and so $\int_{\R^2} h \dx\dy = V$. Then $|\nabla h| = \chi_{[-L,L]^2} H/L$, and $\int_{\R^2} |\nabla h|^2 \dx\dy = 4L^2 (H/L)^2 = 4H^2$. To estimate the elastic energy we observe that for $(x,y,z) \in \Omega_h \cap \{ z \in (0,L) \}$ we have $|\nabla u(x,y,z)|^2 = 2(1-z/L)^2 + (x^2+y^2)/L^2$, and $\nabla u = 0$ otherwise. Since $x^2 + y^2 \le 2L^2$, we have that $\int_{\Omega_h} |\nabla u|^2 \dxyz \le 4|\Omega_h \cap \{ z < L \}| \le 4 (2L)^2 L = 8L^3$, and finally $F(V) \le 4H^2 + 8L^3 = (9/4 + 8)V^{6/7}$. 
 
 \medskip
 It remains to prove the lower bound. First we describe the notation. For $(x_0,y_0) \in \R^2$ and $l > 0$ we define the square $S_l(x_0,y_0) := [x_0,x_0+l) \times [y_0,y_0+l)$. Let $\Phi$ be a function, which for a given square $S_l(x_0,y_0)$ counts on what portion of slices the function $h$ is larger than $l$:
 \begin{equation*}
  \Phi(S_l(x_0,y_0)) := \left|\left\{ x \in [x_0,x_0+l) : h(x,y) \ge l \textrm{ for all } y \in [y_0,y_0+l)\right\}\right| / l.
 \end{equation*}
 We observe that this definition makes sense since $h \in H^1(\R^2)$, and so there exists a representative which is defined everywhere on 
 {{almost every}}
 slice. Let $\eps > 0$ be such that $\int_{\{ h \ge \eps \}} h \dxy \ge V/2$. Given $\eps$, we assume that the following Calder\'on-Zygmund type lemma holds:
 \begin{Lemma}\label{lm:cover}
 There exists a collection of disjoint squares  $\{ S_{l_n}(x_n,y_n) \}$ such that their union covers the set $\{ h \ge \eps \}$, and each square $S_{l_n}(x_n,y_n)$ from the collection satisfies:
 \begin{itemize}
  \item $\Phi(S_{l_n}(x_n,y_n)) \le 1/2$,
  \item there exists a point $(x_n',y_n')$ such that $S_{l_n/2}(x_n',y_n') \subset S_{l_n}(x_n,y_n)$ and $\Phi(S_{l_n/2}(x_n',y_n')) \ge 1/2$.
 \end{itemize}
 \end{Lemma} 

 We postpone the proof of the lemma and first show how the lemma implies the lower bound. Let $S_{l_n}(x_n,y_n)$ be one of the squares obtained in Lemma~\ref{lm:cover}. We denote
 \begin{align*}
S_n &:= S_{l_n}(x_n,y_n), \qquad V_n := \int_{S_n} h \dx\dy ,\qquad S_n' := S_{l_n/2}(x_n',y_n'), 
\\ X_n' &:= \{ x \in [x_n',x_n'+l_n/2) : h(x,y) \ge l_n/2 \textrm{ for all } y \in [y_n',y_n'+l_n/2) \}.
\end{align*} 
Observe that
 \begin{multline}\label{56}
  \int_{S_n\times(0,\infty)\cap \Omega_h} |\nabla u|^2 \dxyz \ge \int_{S_n'\times(0,\infty)\cap \Omega_h} |\nabla u|^2 \dxyz 
\\ \ge \int_{X_n'} \lt( \int_{y_n'}^{y_n'+l_n/2} \int_0^{l_n/2} |\nabla u|^2 \dy \dz \rt) \dx \gtrsim |X_n'| l_n^2 \overset{\Phi(S_n') \ge 1/2}{\gtrsim} l_n^3,
 \end{multline}
 where we used that for $x \in X_n'$ the whole square $\{x\} \times [y_n',y_n'+l_n/2) \times [0,l_n/2) \subset \Omega_h$, and we used the one-dimensional argument to get $\int_{y_n'}^{y_n'+l_n/2} \int_0^{l_n/2} |\nabla u|^2 \dy \dz \gtrsim l_n^2$.

 We consider two cases: $V_n < 2l_n^3$ and $V_n \ge 2l_n^3$. If $V_n < 2l_n^3$, then by~\eqref{56} we have $V_n \lesssim \int_{S_n \times (0,\infty) \cap \Omh} |\nabla u|^2 \dxyz$. 

 Now {{suppose that}}
 $V_n \ge 2l_n^3$. Since $\Phi(S_n) \le 1/2$, we can find $x' \in [x_n,x_n+l_n)$ such that $h(x',y') \le l_n$ for some $y' \in [y_n,y_n+l_n)$ and $\int_{y_n}^{y_n+l_n} |\partial_y h(x',y)|^2 \dy \le 2l_n^{-1}\int_{S_n} |\nabla h|^2 \dx\dy$. 
 Then for any $y \in [y_n,y_n+l_n)$, H\"older's inequality implies
 \begin{multline}
  |h(x',y) - h(x',y')| \le \int_{y}^{y'} |\partial_y h(x',\hat y)| \ud \hat y\\
 \le l_n^{1/2} \left( 2l_n^{-1} \int_{S_n} |\nabla h|^2 \dx\dy \right)^{1/2} = \left( 2\int_{S_n} |\nabla h|^2 \dx\dy \right)^{1/2}.
 \end{multline}
 Since $h(x',y') \le l_n$, we get that $\max_{y \in [y_n,y_n+l_n)} h(x',y) \le l_n + \left( 2\int_{S_n} |\nabla h|^2 \dx\dy \right)^{1/2}$ and that  $\int_{y_n}^{y_n+l_n} h(x',y) \dy \le l_n^2 + l_n\left( 2\int_{S_n} |\nabla h|^2 \dx\dy\right)^{1/2}$. Finally, another application of H\"older's inequality shows that for every $x \in [x_n,x_n+l_n)$:
 \begin{align*}
  \int_{y_n}^{y_n+l_n} h(x,y) \dy &\le \int_{y_n}^{y_n+l_n} h(x',y) \dy + \left| \int_{y_n}^{y_n+l_n} h(x,y) \dy - \int_{y_n}^{y_n+l_n} h(x',y) \dy \right| 
\\ &\le l_n^2 + l_n\left( 2\int_{S_n} |\nabla h|^2 \dx\dy \right)^{1/2} + \int_{S_n} |\nabla h| \dx\dy\\
& \le l_n^2 + 3l_n \left( \int_{S_n} |\nabla h|^2 \dx\dy \right)^{1/2},
 \end{align*}
{{and}} thus
 \begin{equation*}
  V_n = \int_{S_n} h{{\dx\dy}} = \int_{x_n}^{x_n+l_n} \left( \int_{y_n}^{y_n+l_n} h(x,y) \dy \right) \dx \le l_n^3 + 3l_n^2 \left( \int_{S_n} |\nabla h|^2 \dx\dy \right)^{1/2}.
 \end{equation*}
 Since $V_n \ge 2l_n^3$, we have that $V_n^2/(36l_n^4) \le \int_{S_n} |\nabla h|^2 \dx\dy$. Then \eqref{56} and Young's inequality imply
 \begin{equation}\label{512}
  \int_{S_n\times(0,\infty) \cap \Omega_h } |\nabla u|^2 \dxyz + \int_{S_n} |\nabla h|^2 \dx\dy \gtrsim l_n^3 + V_n^2 / (36l_n^4 ) \gtrsim V_n^{6/7}.
 \end{equation}
{{Summarizing, we get }}that 
 \begin{equation*}
  \int_{S_n\times(0,\infty) \cap \Omega_h } |\nabla u|^2 \dxyz + \int_{S_n} |\nabla h|^2 \dx\dy \gtrsim \min{{\{}}V_n^{6/7},V_n{{\}}}.
 \end{equation*}
 Since $\bigcup_{n} S_{l_n}(x_n,y_n) \supset \{ h \ge \eps \}$ and $\int_{ \{ h \ge \eps \}} h \ge V/2$, we have that 
 \begin{equation*}
  \sum_{n} \int_{S_{l_n}(x_n,y_n)} h \dx\dy = \sum_{n} V_n \ge V/2.
 \end{equation*}
 Hence, summing~\eqref{512} over all the squares $S_{l_n}(x_n,y_n)$ and using the concavity of $f(t) := \min{{\{}}t^{6/7},t{{\}}}$ {{yields}}
 \begin{multline*}
  {{F_s^{3D}}}(V) \ge \sum_{n} \left( \int_{S_n\times(0,\infty) \cap \Omega_h } |\nabla u|^2 \dxyz + \int_{S_n} |\nabla h|^2 \dx\dy \right) \\
\gtrsim \sum_{n} \min{{\{}}V_n^{6/7},V_n{{\}}} \gtrsim \min{{\{}}V^{6/7},V{{\}}}.
 \end{multline*} 

\end{proof}

 \begin{proof}[Proof of Lemma~\ref{lm:cover}:]
  Since $h$ is Lipschitz, we see that $\delta := \textrm{dist}( \{ h \ge \eps \}, \{ h \le \eps/2 \}) > 0$. Let $l' := (2V)^{1/3}$, and let $S_{l'}(x_n^{(0)},y_n^{(0)})$ be a (finite) disjoint covering of the set $\{ h \ge \eps \}$. The choice of $l'$ implies that $\Phi(S_{l'}(x_n^{(0)},y_n^{(0)})) \le 1/2$ for every $n$. Let us now fix one square $ S_{l'}(x_n^{(0)},y_n^{(0)})$. We 
  {{aim}} to construct a finite collection {{$\mathcal{N}$}} of disjoint squares with the following properties:
\begin{itemize}
 \item all squares in $\mathcal{N}$ satisfy the two conditions given in Lemma~\ref{lm:cover},
 \item $S_{l'}(x_n^{(0)},y_n^{(0)}) \cap \{ h \ge \eps \} \subset \bigcup_{S \in \mathcal{N}} S$. 
\end{itemize}

We observe that if we construct such $\mathcal{N}$ for each initial square $S_{l'}(x_n^{(0)},y_n^{(0)})$, taking {{the}} union of all those $\mathcal{N}$ will give a collection of squares which 
{{satisfies}} all the required conditions of Lemma~\ref{lm:cover}.

We  now describe the iterative construction. We define $\mathcal{M}_0 := \{ S_{l'}(x_n^{(0)},y_n^{(0)}) \}$ and $\mathcal{N}_k := \emptyset$. Assume that $\mathcal{M}_k$ is constructed, and let $S_l(x,y) \in \mathcal{M}_k$. We divide $S_l(x,y)$ into four disjoint squares $S_1,S_2,S_3,S_4$ with sidelength $l/2$, and consider two cases. If $\Phi(S_i) \le 1/2$ for all $i=1,2,3,4$, then we add to $\mathcal{M}_{k+1}$ all those $S_i$ which satisfy $S_i \cap \{ h \ge \eps \} \neq \emptyset$. Otherwise (i.e., if $\Phi(S_i) > 1/2$ for some $i$), we add $S_l(x,y)$ into $\mathcal{N}_k$. 

First, we observe that at any step $k_0$ of the procedure we have 
\begin{equation}
 S_{l'}(x_n^{(0)},y_n^{(0)}) \cap \{ h \ge \eps \} \subset \bigcup_{S \in \mathcal{M}_{k_0} \cup \bigcup_{k < k_0} \mathcal{N}_k} S.
\end{equation}
We also see that any square $S \in \mathcal{N}_k$ satisfies the two conditions given in Lemma~\ref{lm:cover}, and so we are done with the construction of $\mathcal{N} := \bigcup N_k$ provided we show that $\mathcal{M}_k = \emptyset$ for large enough $k$. Let $k$ be such that $l' 2^{-k} \le \min {{\{}}\delta/\sqrt{2},\eps/2 {{\}}}$, and let us assume that $S \in \mathcal{M}_k$. From the way $\mathcal{M}_k$ was constructed we see that $\Phi(S) \le 1/2$, the length of the side of $S$ is $l'2^{-k}$, and $S \cap \{ h \ge \eps \} \neq \emptyset$. 

On the other hand, since $l'2^{-k} \le \delta/\sqrt{2}$, we see that $S \subset \{ h \ge \eps/2 \}$. Then, $l'2^{-k} \le \eps/2$ implies $h \ge l'2^{-k}$ in $S$, and so $\Phi(S) = 1$, which contradicts $\Phi(S)\le 1/2$. 

We showed that $\mathcal{M}_k$ is empty for large enough $k$, and so the $\mathcal{N}$ we just constructed satisfies all the requirements of Lemma~\ref{lm:cover}. This finishes the proof of 
{{Lemma~\ref{lm:cover}, and thus of Theorem~\ref{th:scaling3D}.} }

 \end{proof}

Using the same line of proof, we can also show the scaling law for the large slope approximation, i.e., for the functional
\begin{eqnarray*}
 {F}_\ell^{3D}(V):=\inf_{(u,h)} \left\{\int_{\Omega_h} |\nabla u|^2 \dxyz + \int_{\R^2} |\nabla h| \dx\dy : {{h\in BV(\R^2),\,h \in W^{1,\infty}(\R^2),\ u\in H^1(\Omega_h)}}\right.\\
\left. \int_{\R^2} h \dx\dy=V,  \, u(x,y,0)=(x,y) \right\},
\end{eqnarray*}
\begin{Theorem}
 There exists a positive constant $c$ such that for every $V > 0$
 \begin{equation}
  c \min {{\{}}V,V^{3/4}{{\}}} \le {{F_\ell^{3D}}}(V) \le\frac{1}{c} \min{{\{}}V,V^{3/4}{{\}}}.
 \end{equation} 
\end{Theorem}
Notice that in contrast to \cite[Proposition 5.1]{GolZwick}, we do not need here any hypothesis on $V$.
%

}


\section*{Acknowledgments}
The authors would like to thank F. Otto for suggesting the problem and for interesting discussions. M. Goldman was funded by the Von Humboldt postdoctoral fellowship. 
The authors gratefully acknowledge the hospitality of the MPI MiS Leipzig (B. Zwicknagl) and of the SFB 1060  at the University of Bonn (M. Goldman), where part of this work has been carried out.
\bibliography{islands}
\end{document}